\RequirePackage{fix-cm}
\documentclass[twocolumn]{svjour3}          % twocolumn
\smartqed  % flush right qed marks, e.g. at end of proof
\usepackage{graphicx}
\usepackage{graphicx}

\usepackage{tikz-cd}
\usepackage{tikz}

\usetikzlibrary{decorations.markings}
\usetikzlibrary{shapes.geometric, arrows}
\usepackage{tikzscale}
\usepackage{filecontents}
\usepackage{wrapfig}
\usepackage{tcolorbox}
\usepackage{color}
\definecolor{blue1}{rgb}{0.0, 0.0, 1.0}
\definecolor{gray}{rgb}{0.9,0.9,0.9}
\definecolor{gray1}{rgb}{0.7,0.7,0.7}
\definecolor{gray2}{rgb}{0.8,0.8,0.8}
\definecolor{magenta}{rgb}{1.0, 0.0, 1.0}

\usepackage{hyperref}
\hypersetup{ colorlinks=true, urlcolor  = blue, linkcolor = blue,
citecolor = blue1,}

\usepackage{float}
\usepackage{amsmath,mathtools,amssymb,amsfonts}
\usepackage[scaled=.95]{helvet}

 \journalname{X - Journal}
\begin{document}

\title{Weak-noise-induced transitions with inhibition and modulation of neural oscillations}


\author{Marius E. Yamakou \and J\"{u}rgen Jost
}


\institute{M. E. Yamakou $\cdot$ J. Jost \at Max-Planck-Institut
f\"{u}r Mathematik in den Naturwissenschaften, Inselstr. 22, 04103 Leipzig, Germany\\
M. E. Yamakou $\cdot$ J. Jost \at Fakult\"{a}t f\"{u}r Mathematik und Informatik, Augustusplatz 10,
04109 Leipzig, Germany\\
J. Jost \at Santa Fe Institute for the Sciences of Complexity,\\ NM 87501, Santa Fe, USA\\
             \email{yamakou@mis.mpg.de}  \emph{Corresponding author}\\
             \email{jost@mis.mpg.de}           %  \\
            %of F. Author  %  if needed
}

\date{Received: date / Accepted: date}

\maketitle

\begin{abstract}
We analyze the effect of weak-noise-induced transitions on the 
dynamics of the FitzHugh-Nagumo neuron model in a bistable state consisting of a stable 
fixed point and a stable unforced limit cycle. 
Bifurcation and slow-fast analysis give conditions on the parameter 
space for the establishment of this bi-stability. 
In the parametric zone of bi-stability, weak-noise amplitudes may strongly 
inhibit the neuron's spiking activity. Surprisingly, increasing the noise strength leads
to a minimum in the spiking activity, after which the activity starts to increase monotonically with 
increase in noise strength. We investigate this inhibition and 
modulation of neural oscillations by weak-noise amplitudes by looking at 
the variation of the mean number of spikes per unit time with the noise intensity.
We show that this phenomenon always occurs when the initial conditions
lie in the basin of attraction of the stable limit cycle.
For initial conditions in the basin of attraction of the stable fixed point, 
the phenomenon however disappears, unless the time-scale separation parameter 
of the model is bounded within some interval. We provide a theoretical 
explanation  of this phenomenon in terms of the stochastic sensitivity 
functions of the attractors and their minimum Mahalanobis distances from the 
separatrix isolating the basins of attraction. 
\end{abstract}

%%%%%%%%%%%%%%%%%%%%%%%%%%%%%%%%%%%%%%%%%%%%%%%%%%%%%%%%%%%%%%%%%%%%%%%%%%%%%%%%%%%%%%%%%%%%%%%%%%%%%%%%%%%%%%%%%%%%%%%%%
\section{Introduction}\label{sect1}
%%%%%%%%%%%%%%%%%%%%%%%%%%%%%%%%%%%%%%%%%%%%%%%%%%%%%%%%%%%%%%%%%%%%%%%%%%%%%%%%%%%%%%%%%%%%%%%%%%%%%%%%%%%%%%%%%%%%%%%%%
\noindent
Fixed points, periodic, quasi-periodic or chaotic orbits 
are typical solutions of deterministic nonlinear dynamical
systems. Multi-stability is common, as  a dynamical
system typically  possesses two or more mutually exclusive stable
solutions (attractors). For a given set of the parameters,
coexistent stable states represented by different (or identical) types
of attractors in the phase space of the system are by topological
necessity separated by some unstable states. In neurodynamics for
example, spiking neurons may possess coexistent quiescent (fixed
point) and tonic spiking states (limit cycle) \cite{Paydarfar}, or
distinct periodic and chaotic spiking states \cite{G. Cymbalyuk}.
A given state can be reached if the system starts from a set of
initial conditions within the state's basin of attraction.
Otherwise, an external perturbation can be used to switch the
system from one stable attractor to another. When noise is introduced into
the system, random trajectories can visit different stable states 
of the system by jumping over the unstable ones. 

Important and challenging problems in multi-stable systems are to
find the residence times of random trajectories in each stable state and
its statistics, and the critical value of the noise amplitude and
control parameters at which noise-induced jumping becomes
significant. The analytical treatment of such problems based on the
Fokker-Planck Equation (FPE) becomes complicated for $n$-dimensional
dynamical systems, $n\geq2$, and therefore, various approximations
were developed and are now commonly used \cite{N. G. Van
Kampen,M.I Freidlin}. 

The quasi-potential method gives exponential
asymptotics for the stationary probability density. In the
vicinity of the deterministic attractor, the first approximation
of the quasi-potential is a quadratic form \cite{G. Mil'shtein},
leading to a Gaussian approximation of the stationary probability
density of the FPE. The corresponding covariance matrix characterizes the 
stochastic sensitivity of the deterministic attractor: its
eigenvalues and eigenvectors define the geometry of bundles of
stochastic trajectories around the deterministic attractors.
The Gaussian distribution centered on an attractor can be viewed as a
confidence ellipsoid, while a minimal distance from this ellipsoid
to the boundary separating the basins of attraction is proportional to
the escape probability \cite{I. Bashkirtseva}. The appropriate
measure for this distance is the so-called Mahalanobis distance
\cite{P. C. Mahalanobis}, the distance from a point to a
distribution.

The residence time of random trajectories in a basin of attraction
depends on two factors. The first factor is the  geometry 
of the  basin of attraction, e.g., the larger the
distance is between an attractor and the separatrix isolating its
basin of attraction, the longer is the residence time of phase
trajectories in the basin. Second, the 
attractors are sensitive to random perturbations: the higher the stochastic
sensitivity function (SSF) is, the higher is the probability to
escape from the basin of attraction, and thus the shorter is the
residence time \cite{M.I Freidlin}. Therefore, considering only
the geometrical arrangement of stable attractors and the
separatrix (the Euclidean distance between them) might not be
sufficient for a theoretical explanation of a stochastic
phenomenon, like the one  to be investigated in the present work, and so the
sensitivity of the attractors to random perturbations must also be
taken into account. The Mahalanobis distance, which combines
geometrical and stochastic sensitivity aspects of the dynamics,
allows for a proper theoretical explanation of the transitions
between attractors.

The effects of
noise in neurobiological dynamical systems have been intensively investigated, for both single
neurons and neural networks. Some of the most
studied noise-induced phenomena are: stochastic resonance (SR)
\cite{Lindner,Longtin,Collins}, coherence resonance (CR)
\cite{Pikovsky}, and noise-induced synchronization
\cite{Kim}. During SR, the neuron's spiking activity becomes
more closely correlated with a sub-threshold periodic input current
in the presence of an optimal level of noise. In $1997$, Pikovsky
and Kurths showed that CR is basically SR in the absence of a
periodic input current. During CR, noise can activate the
neuron by producing a sequence of pulses which can achieve a
maximal degree of coherence for an optimal noise amplitude if the system is
in the neighborhood of its Andronov-Hopf bifurcation. We notice in
these phenomena that noise has a facilitatory effect and leads only to
increased responses.

More recently, it was
discovered both experimentally \cite{Paydarfar} and theoretically Gutkin et \textit{al.}
\cite{Gutkin1,Gutkin2} (see also \cite{Jost}) that noise can also  turn off repetitive
neuronal activity . 
\cite{Gutkin1,Gutkin2,Jost} used the Hodgkin-Huxley equations in bistable regime
(fixed point and limit cycle)
with a mean input current consisting of both a deterministic 
and random input component, to computationally confirm the
inhibitory and modulation effects of Gaussian noise on the neuron's 
spiking activity. 
They found that there is a tuning effect of noise that has the
opposite character to SR and CR, which they termed inverse
stochastic resonance (ISR). Very recently (August 2016),
the first experimental confirmation of ISR and its plausible functions in 
local optimal information transfer was reported in \cite{Buchin}, 
where the Purkinje cells that play a central role in the cerebellum are 
used for the experiment.
During ISR, weak-noise amplitudes may strongly inhibit the spiking 
activity down to a minimum level (thereby decreasing the mean 
number of spikes to a minimum value), after which the activity 
starts and continuously increases with increasing noise amplitude
(thereby monotonically increasing the mean number of spikes with 
increasing noise amplitude). In \cite{Gutkin1,Gutkin2}, it was shown 
that ISR  occurred and persisted regardless of which  basin of attraction the 
initial conditions are chosen from, \textit{provided} the deterministic 
input current component is above its Andronov-Hopf bifurcation value. 

In the present work, we investigate ISR in a theoretical neuron model in
the absence of a deterministic input current component. We consider
the model with only a random input component, which is
in a state of bi-stability consisting of a stable fixed point and a
stable unforced limit cycle. We show that ISR occurs as well in
this case and greatly depends not only on the location of the initial conditions,
but also on the time-scale separation
parameter of the model. More precisely, we show that ISR \textit{always} occurs when the 
initial conditions are chosen from the basin of attraction of the stable
limit cycle. When the initial conditions are in the basin of
attraction of the stable fixed point, we show that ISR disappears, 
except interestingly when
the time-scale separation parameter of the model lies within a certain interval. 
A theoretical explanation of this phenomenon is given in terms of 
the SSFs of the stable attractors and their 
Mahalanobis distances from the separatrix.

This paper is organized as follows: In Sect.\ref{sect2}, we
present the theoretical neuron model used to analyze ISR.
In Sect.\ref{sect3}, we make explicit deterministic bifurcation
and slow-fast analyses of the model equation and show how
bi-stability consisting a stable fixed point and a stable unforced
limit cycle establishes itself. In Sect.\ref{sect4}, we make a
stochastic sensitivity analysis of the stable attractors. In
Sect.\ref{sect5}, we investigate ISR through numerical simulations
and provide a theoretical explanation of the phenomenon using the
results in Sect.\ref{sect4}. In
Sect.\ref{sect6}, we have concluding remarks.

%%%%%%%%%%%%%%%%%%%%%%%%%%%%%%%%%%%%%%%%%%%%%%%%%%%%%%%%%%%%%%%%%%%%%%%%%%%%%%%%%%%%%%%%%%%%%%%%%%%%%%%%%%
\section{Model equation and phenomenon}\label{sect2}
%%%%%%%%%%%%%%%%%%%%%%%%%%%%%%%%%%%%%%%%%%%%%%%%%%%%%%%%%%%%%%%%%%%%%%%%%%%%%%%%%%%%%%%%%%%%%%%%%%%%%%%%%%%
\noindent
In this paper, we consider a stochastic perturbation of a version of the Fitzhugh-Nagumo (FHN) neuron model 
\cite{FitzHugh}. We consider the resulting stochastic differential equation both in
the slow time-scale $\tau$ (Eq.\eqref{fn0a}) and  in the fast
time-scale $t$ (Eq.\eqref{fn0b})
\begin{equation}\label{fn0a}
\begin{split}
\left\{\begin{array}{lcl}
dv_{\tau}&=&\frac{1}{\varepsilon}f(v_{\tau},w_{\tau})d\tau+\frac{\sigma}{\sqrt{\varepsilon}}dW_{\tau},\\
dw_{\tau}&=&g(v_{\tau},w_{\tau})d\tau,
\end{array}\right.
\end{split}
\end{equation}
\begin{equation}\label{fn0b}
\begin{split}
\left\{\begin{array}{lcl}
dv_{t}&=&f(v_{t},w_{t})dt+\sigma dW_{t},\\
dw_{t}&=&\varepsilon g(v_{t},w_{t})dt,
\end{array}\right.
\end{split}
\end{equation}
with the deterministic velocity vector fields given by
\begin{equation}\label{2}
\begin{split}
\left\{\begin{array}{lcl}
 f(v,w)&=&v(a-v)(v-1)-w, \\
 g(v,w)&=&bv-cw ,
 \end{array}\right.
\end{split}
\end{equation}
where $(v,w)\in\mathbb{R}^2$ represent the activity of the
action potential $v$ and the recovery current $w$ that restores
the resting state of the model. We have as constant parameters
$b>0$, $c>0$, and $a$ is often confined to the range $0<a<1$, but
the case $a < 0$ will be examined in this work.

We have a singular parameter, $0<\varepsilon:=\tau/t\ll1$, i.e.,
the time-scale separation ratio between the slow time-scale $\tau$
and the fast time-scale $t$. We note that Eq.\eqref{fn0a}
preserves the sense of the dynamics on the trajectories of
Eq.\eqref{fn0b}. In other words, the phase trajectories of both
systems of dynamical equations have exactly the same dynamical
behavior. The only difference is the speed of these trajectories
in the phase space. Because the speeds of the trajectories do not
affect in any way our analysis, we will work on both time-scales.
The slow time-scale equation at some points allows for quicker
conclusions in bifurcation analysis while the fast time-scale
equation has an advantage in numerical simulations as it
avoids the division by the very small parameter $\varepsilon$.

$dW_t$ is standard white noise, the formal derivative of
Brownian motion with mean zero and unit variance, and $\sigma$ is
the amplitude of this noise. The random term in Eq.\eqref{fn0a} is rescaled in Eq.\eqref{fn0b} 
according to the scaling law of Brownian motion. That is,
if $ W_t$ is a standard Brownian motion, then for every 
$\lambda\geq0$, $\lambda^{-1/2} W_{\lambda t}$ is also a standard
Brownian motion, i.e., the two processes have the same distribution \cite{Durrett}.

\begin{figure}%[H]
\begin{center}
\textbf{(a)}\includegraphics[width=8.0cm,height=4.1cm]{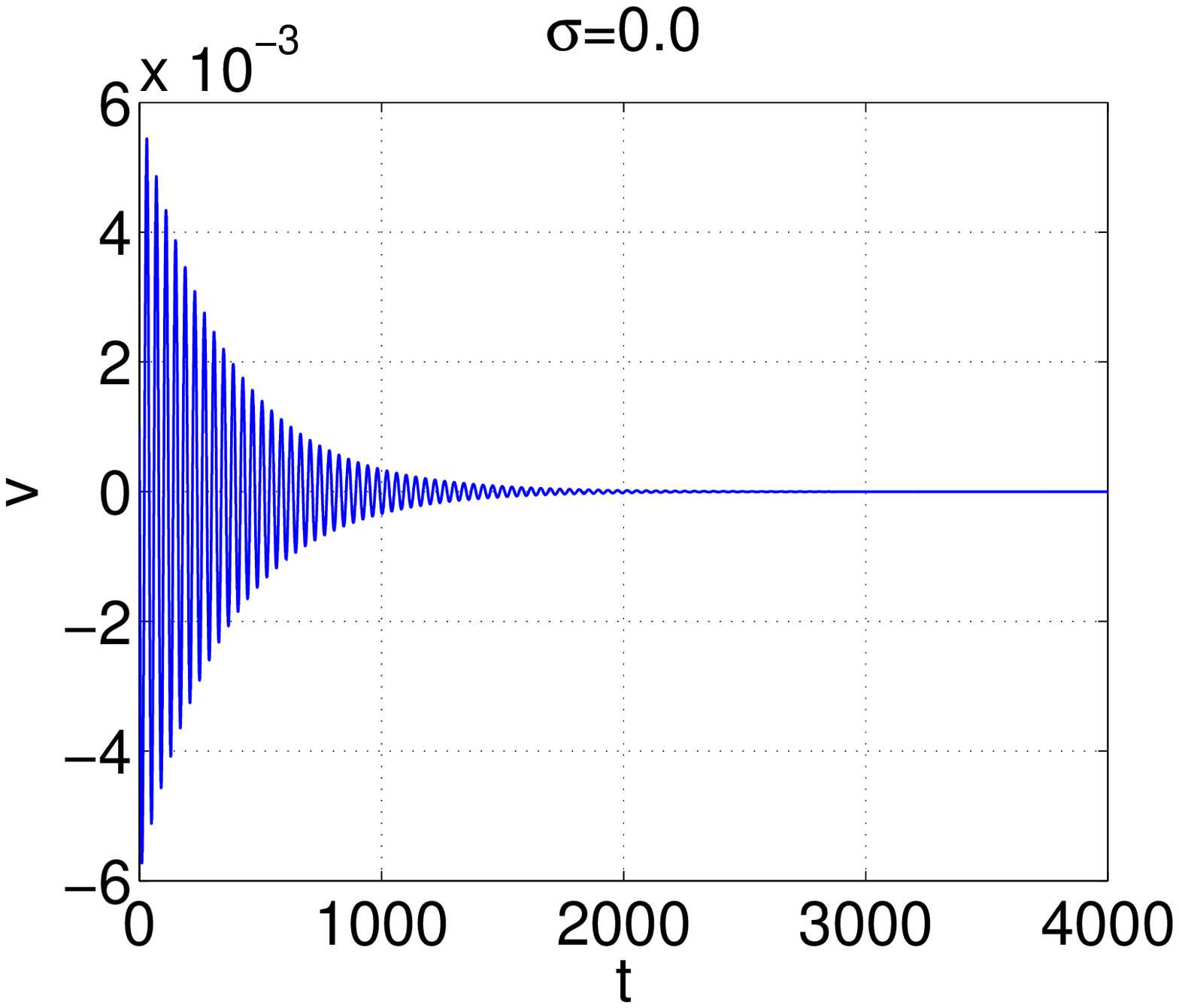}
\textbf{(b)}\includegraphics[width=8.0cm,height=4.1cm]{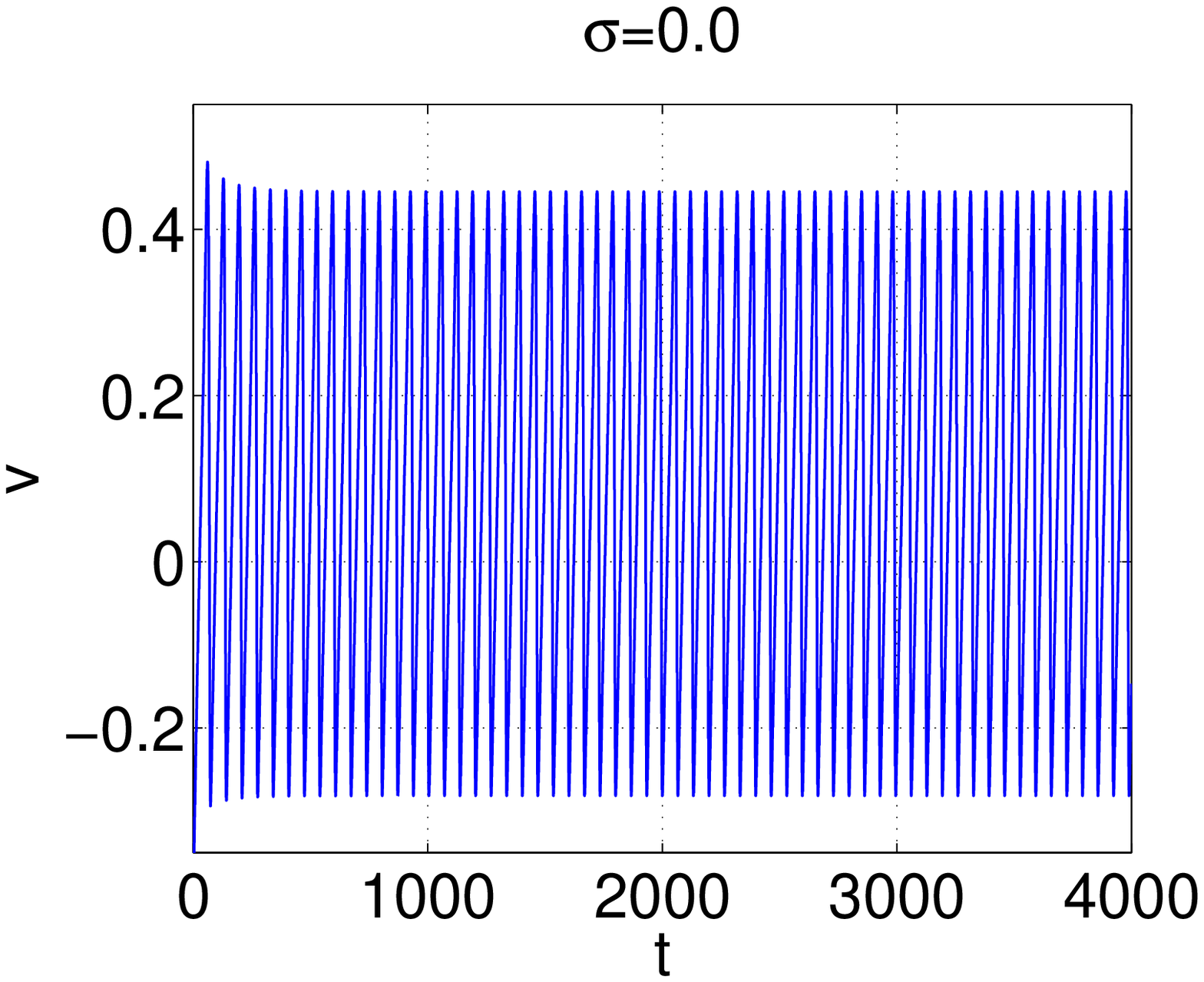}
\textbf{(c)}\includegraphics[width=8.0cm,height=4.1cm]{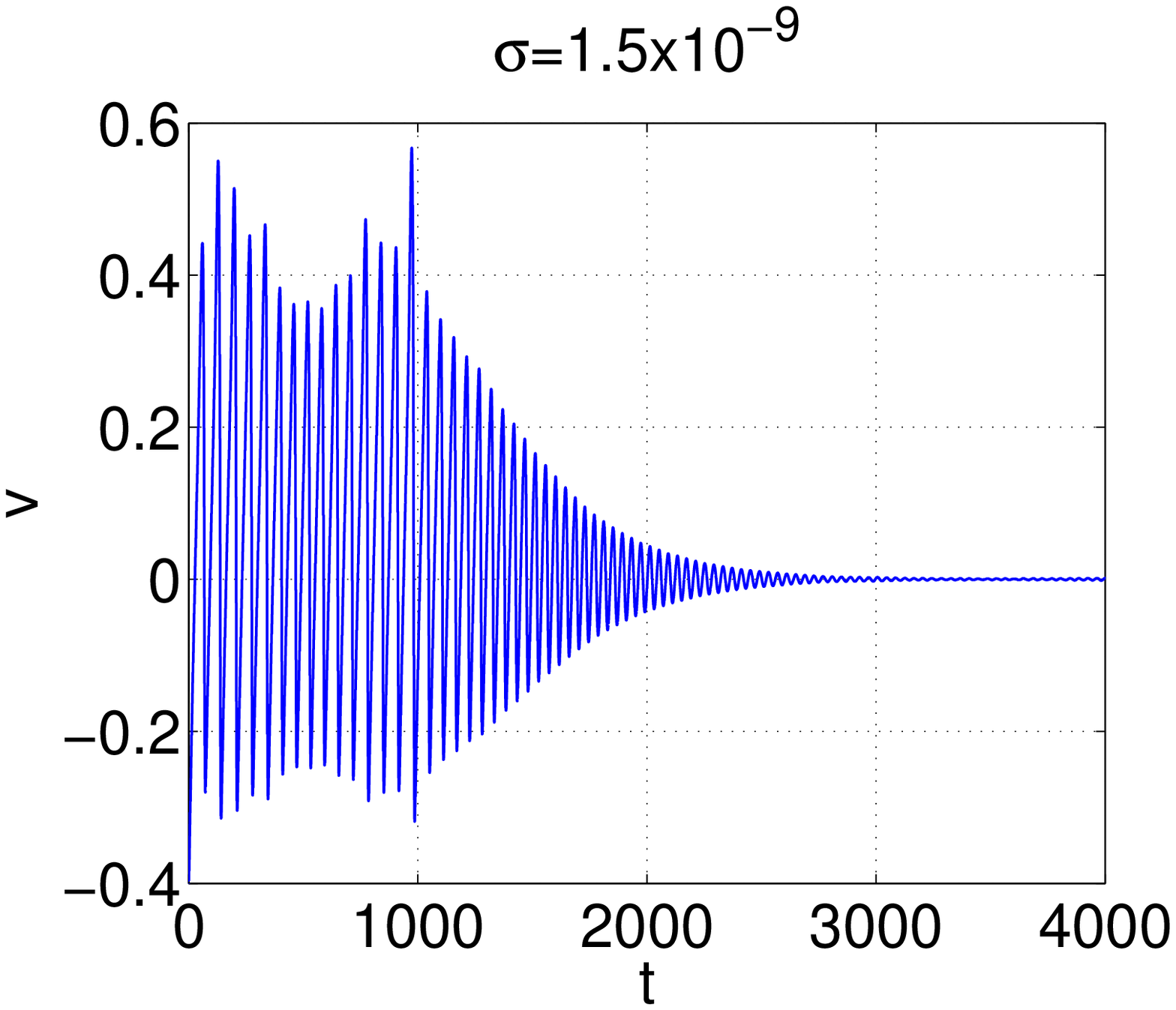}
\textbf{(d)}\includegraphics[width=8.0cm,height=4.1cm]{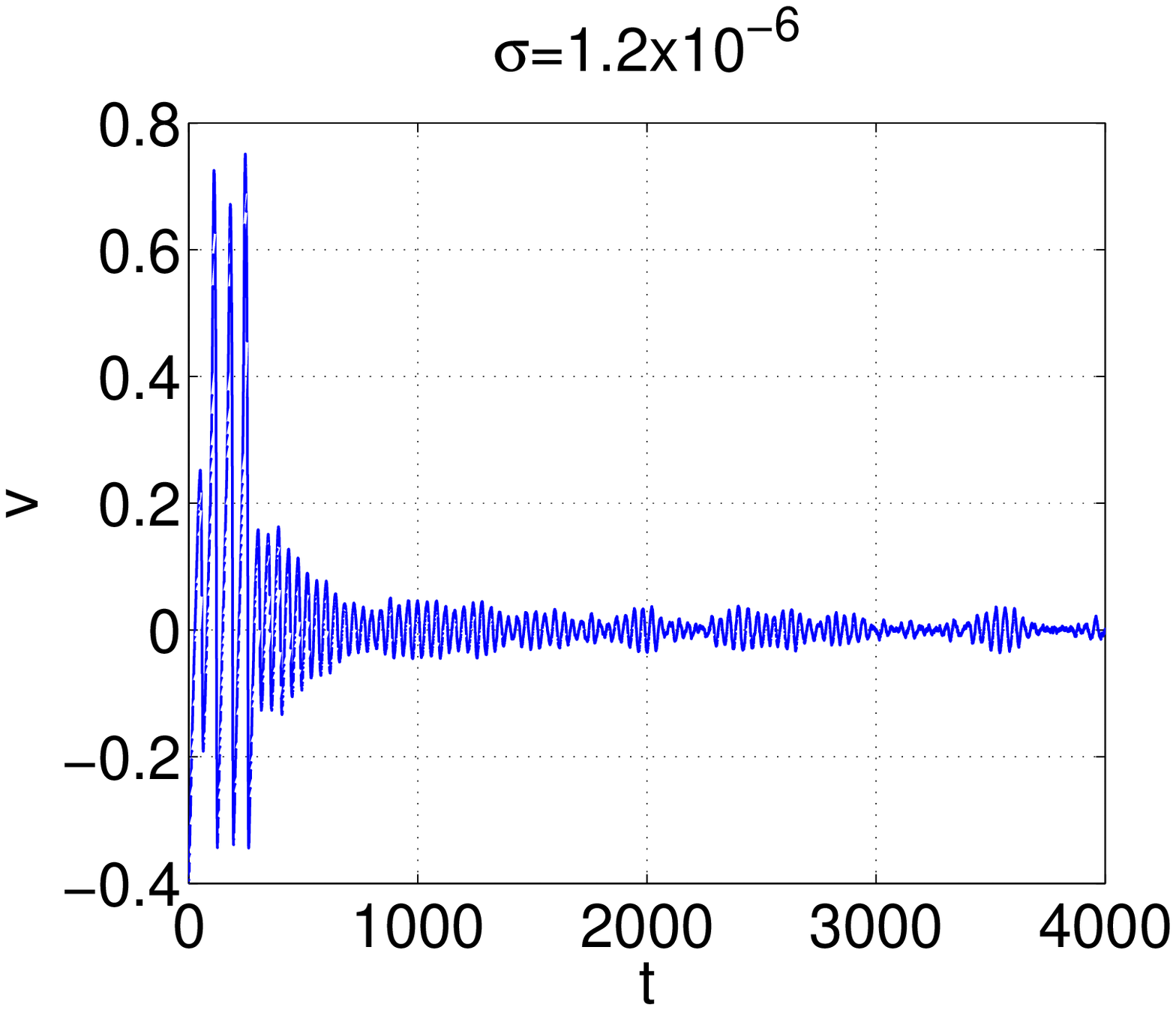}
\textbf{(e)}\includegraphics[width=8.0cm,height=4.1cm]{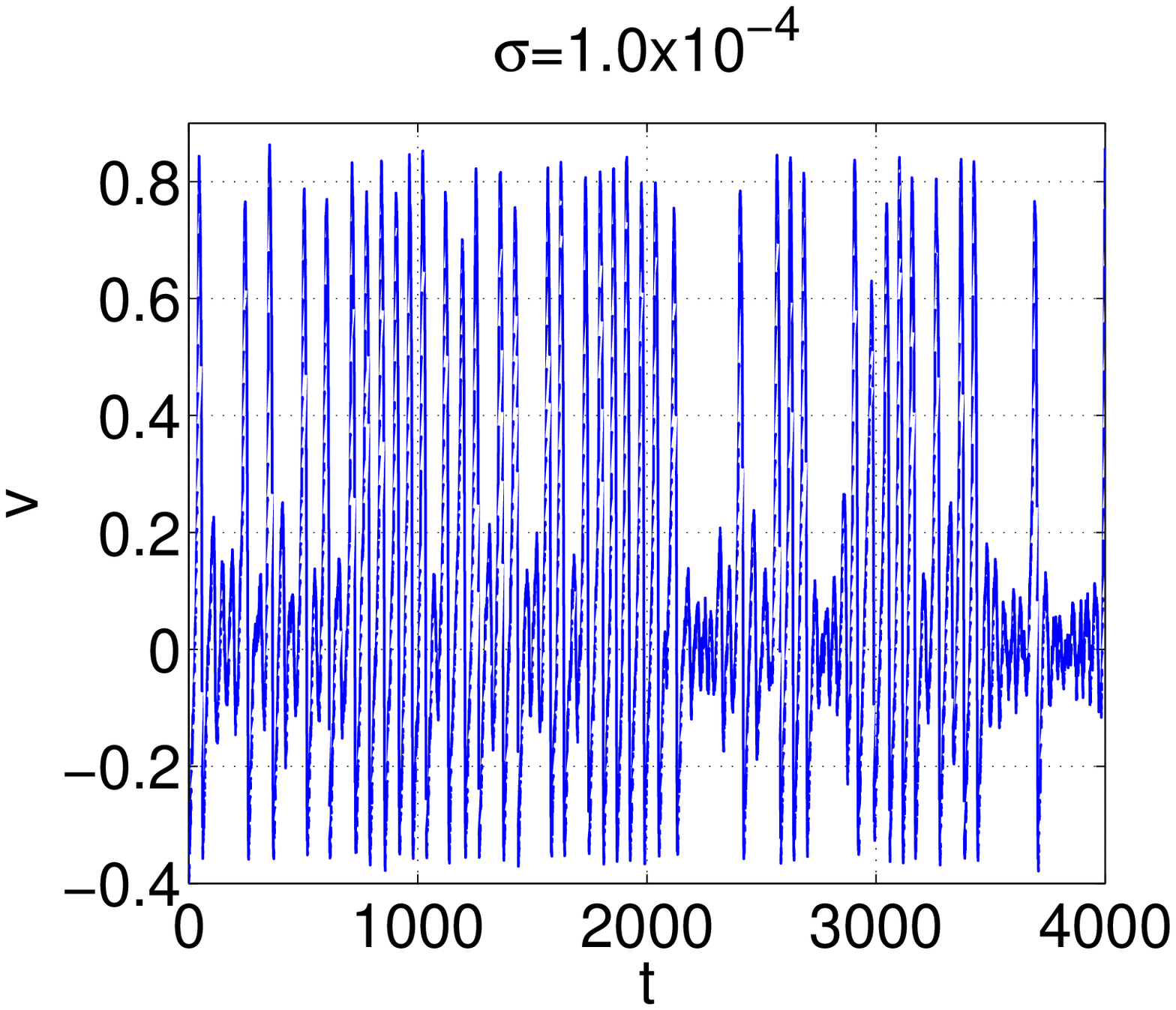}
\caption{Time series of the action potential variable $v$ in Eq.\eqref{fn0b}. \textbf{(a)} and
\textbf{(b)} show the zero-noise dynamics with initial condition in \textbf{(a)} at
$\big(v(0),w(0)\big)=(0.001,0.001)$ and in \textbf{(b)} at $\big(v(0),w(0)\big)=(-0.4,0.2)$.
\textbf{(c)}, \textbf{(d)}, and \textbf{(e)} show the effects of noise of various intensities, $\sigma$, on
the dynamics of $v$ with $\big(v(0),w(0)\big)=(-0.4,0.2)$ in each case.
$a=-0.05, b=1.0, c=2.0, \varepsilon=0.02785$ } \label{fig:chap31}
\end{center}
\end{figure}

Fig.\ref{fig:chap31} shows the time series produced by 
the dynamics of the action potential variable $v$.
In the deterministic case ($\sigma=0$), and for
$a=-0.05$, $b=1.0$, $c=2.0$, and $\varepsilon=0.02785$,
Eq.\eqref{fn0b} can result in  two different dynamics.
In Fig.\ref{fig:chap31}\textbf{a} with initial conditions at
$\big(v(0),w(0)\big)=(0.001,0.001)$, the neuron has only sub-threshold
oscillations with $v$ converging  to zero
and remaining at this value as the time $t$ increases.
In Fig.\ref{fig:chap31}\textbf{b} with initial conditions now at
$\big(v(0),w(0)\big)=(-0.4,0.2)$, the neuron shows self-sustained
supra-threshold oscillations. Thus, the system is bistable.

Fig.\ref{fig:chap31}\textbf{c}-\textbf{e} show a stochastic behavior
($\sigma>0$), with again $a=-0.05$, $b=1.0$, $c=2.0$,
$\varepsilon=0.02785$, and the initial conditions all at
$\big(v(0),w(0)\big)=(-0.4,0.2)$. We
count a spike when $v$ is greater than or equal to the threshold
value $v_{th}=0.25$. In Fig.\ref{fig:chap31}\textbf{c} with a weak-noise
amplitude, i.e., $\sigma=1.5\times10^{-9}$, we have supra-threshold
oscillations for a certain time length with  $21$ spikes after
which $v$ starts to converge to zero and the spiking eventually
stops. In Fig.\ref{fig:chap31}\textbf{d}, when the noise
amplitude is increased (but still relatively weak) to
$\sigma=1.2\times10^{-6}$, we have an even faster inhibition of the spiking with a smaller number of spikes.
In this realization, we have only $3$
spikes. In Fig.\ref{fig:chap31}\textbf{e}, with a stronger noise
amplitude, $\sigma=1.0\times10^{-4}$, the number of spikes
increases again up to $38$. Intuitively, it is surprising that
weak-noise amplitudes inhibit the spiking activity of the neuron
with the occurrence of a minimum in the number of spikes as the noise
amplitude increases even though the initial conditions are exactly
the same as in Fig.\ref{fig:chap31}\textbf{b}. 
We shall investigate in detail the mechanisms behind this phenomenon.
%%%%%%%%%%%%%%%%%%%%%%%%%%%%%%%%%%%%%%%%%%%%%%%%%%%%%%%%%%%%%%%%%%%%%%%%%%%%%%%%%%%%%%%%%%%%%%%%%%%%%%%%%%%
\section{Bifurcation and slow-fast analysis}\label{sect3}
\begin{figure}%[H]
\begin{center}
\textbf{(a)}\includegraphics[width=7.5cm,height=7.0cm]{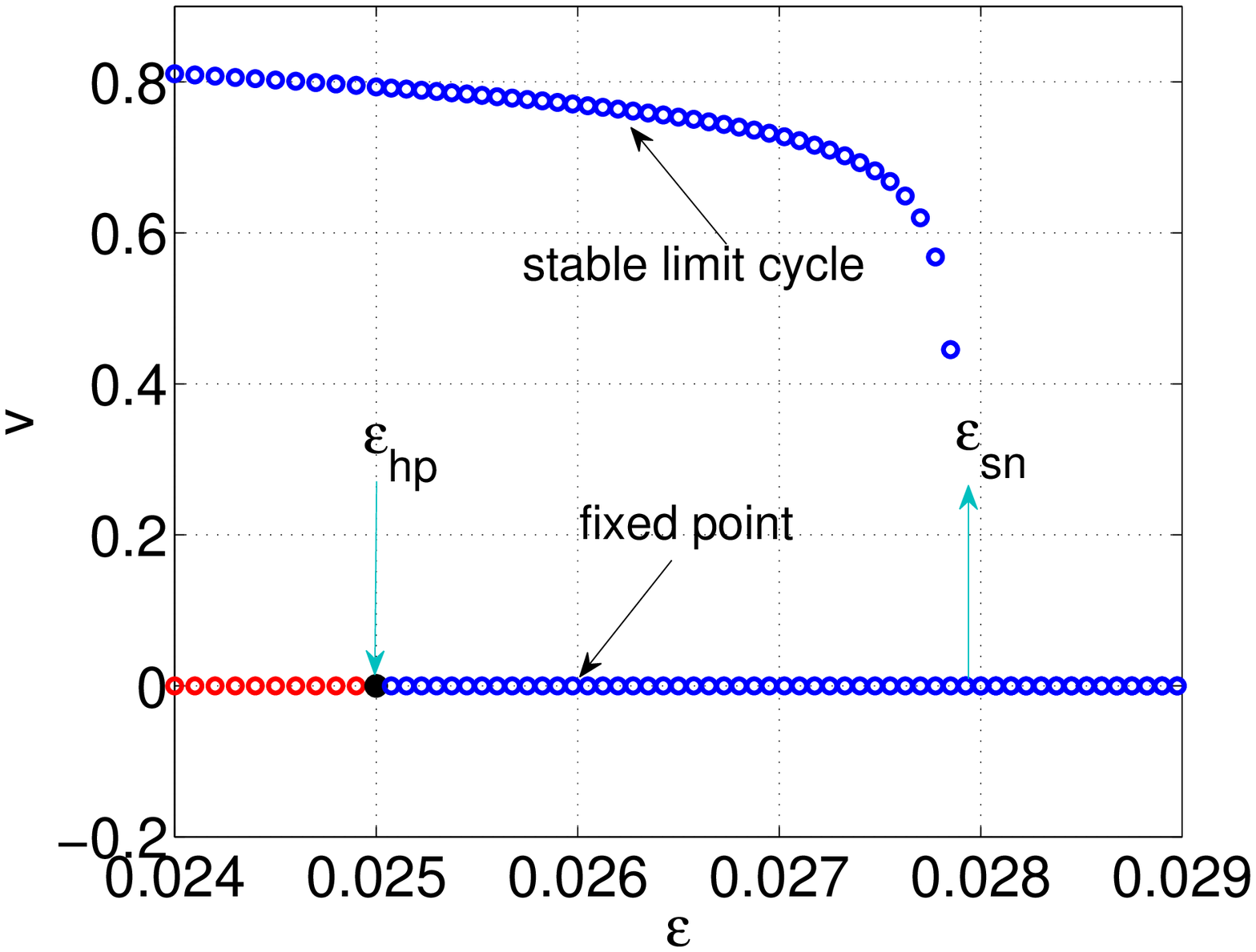}
\textbf{(b)}\includegraphics[width=7.5cm,height=7.0cm]{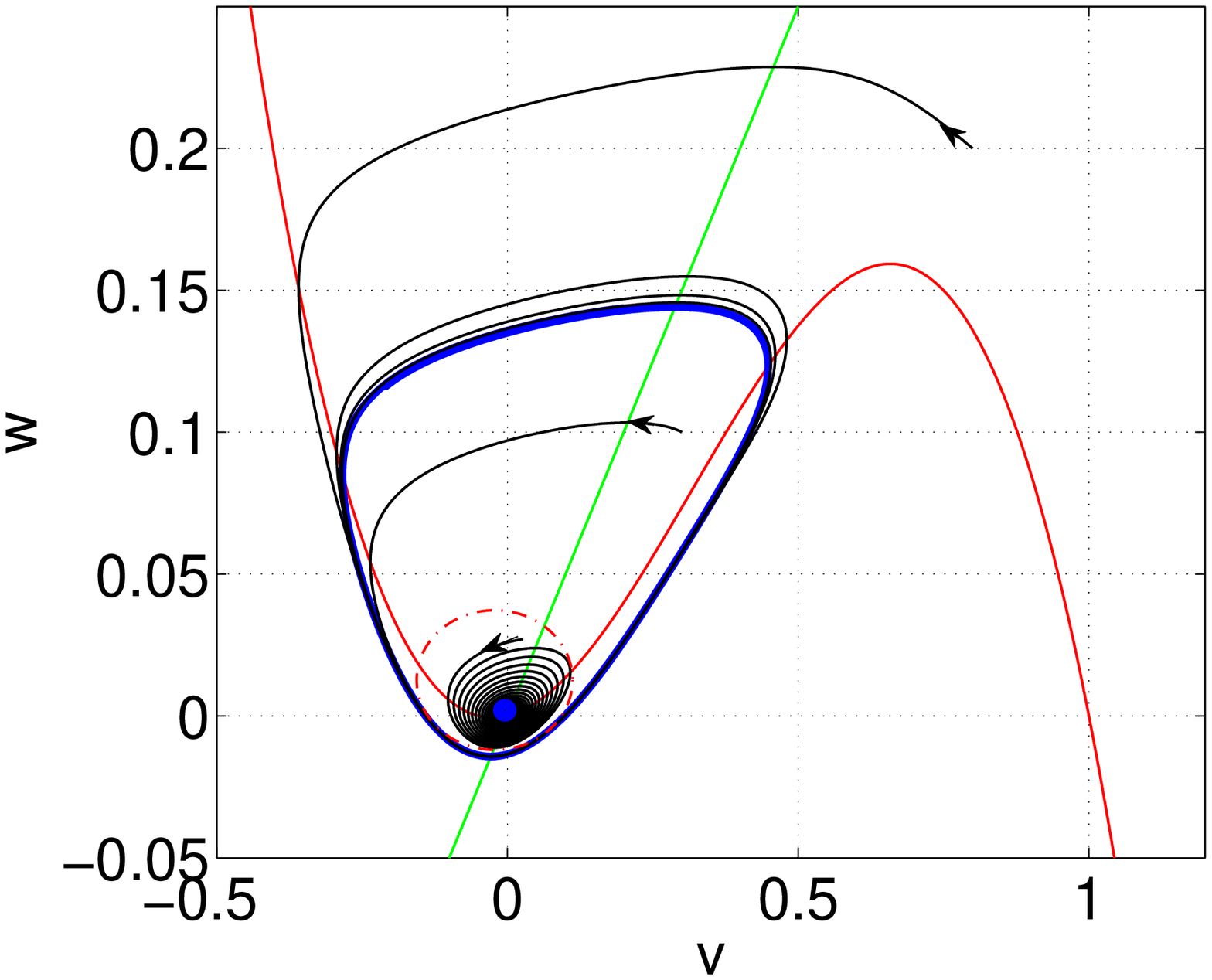}
\caption{(Color online) \textbf{(a)} Bifurcation diagram for Eq.\eqref{fn0b} with a fixed
point at $(v_0,w_0)=(0,0)$ unstable in the singular parameter range 
$0.024<\varepsilon<0.025$ shown by the red points
and a stable limit cycle in this same interval. At
$\varepsilon=\varepsilon_{hp}=0.025$, $(v_0,w_0)$ gain stability
through a sub-critical Andronov-Hopf bifurcation and therefore
co-exist with the stable limit in the interval
$0.025<\varepsilon<0.027865$.  The stable limit cycle undergoes
a saddle-node bifurcation and disappears by shrinking and
eventually colliding with the boundary of the basin of attraction
(the unstable limit cycle not shown) of the stable fixed point at
$\varepsilon=\varepsilon_{sn}=0.027865$. In
$0.027865<\varepsilon\leq0.029$, there is only the stable fixed
point $(v_0,w_0)$ in the entire phase space. \textbf{(b)} Geometry of attractors of
Eq.\eqref{fn0b} with $\varepsilon_{hp}<\varepsilon<\varepsilon_{sn}$. 
The red curve represents the cubic critical
manifold $\mathcal{M}_0$ intersecting the $w$-nullcline (the green
line) at the blue dot corresponding to the fixed point $(v_0,w_0)=(0,0)$, located
to the right of the minimum of $\mathcal{M}_0$ at $v_-=-0.25305$. The blue closed curve
represents the stable limit cycle, the red dotted closed curve the
separatrix (unstable limit cycle), and 3 different trajectories (in black) with arrows at
the initial conditions. Depending on which side of the separatrix
the initial conditions are chosen, solutions converge to either the
stable limit cycle or to the stable fixed point. $a=-0.05, b=1.0,
c=2.0, \varepsilon=0.02785, \sigma=0.0$} \label{fig:chap32}
\end{center}
\end{figure}

%%%%%%%%%%%%%%%%%%%%%%%%%%%%%%%%%%%%%%%%%%%%%%%%%%%%%%%%%%%%%%%%%%%%%%%%%%%%%%%%%%%%%%%%%%%%%%%%%%%%%%%%%%%
\noindent
We now consider the deterministic dynamics corresponding to  Eq.\eqref{fn0a} when $\sigma = 0$ and perform an explicit
bifurcation and slow-fast analysis through which we find the parametric conditions for the
establishment of a bi-stability regime consisting of a stable fixed
point and a stable unforced limit cycle. We note again that our model has no
deterministic input current. At the fixed points $(v_e,w_e)$
(rest states of the neuron), the variables $v(t)$ and $w(t)$ reach
a stationary state, while the set of fixed points is defined
by the intersection of nullclines as
\begin{equation}\label{fn1}
(v_e,w_e):=\Big\{(v,w)\in\mathbb{R}^2:f(v,w)=g(v,w)=0\Big\}.
\end{equation}
From Eq.\eqref{fn1}, we obtain the fixed point equations
\begin{align}\label{fn20}
\begin{split}
\left\{\begin{array}{lcl}
 \frac{b}{c}v&=&-v^3+(a+1)v^2-av,\\
  w&=&\frac{b}{c}v,
 \end{array}\right.
\end{split}
\end{align}
which has the solutions for $v$ as
\begin{equation}\label{fn21}
v_0=0 \text{ and } v_{1,2}=\frac{a+1}{2}\pm \sqrt{\frac{(a-1)^2}{4} -\frac{b}{c}},
\end{equation}
where the solutions $v_{1,2}$ exist if
\begin{equation}\label{fn22}
\frac{(a-1)^2}{4}\ge \frac{b}{c}.
\end{equation}
We always assume $b,c >0$.
When $0<a<1$, we then have
\begin{equation}\label{fn22a}
v_0 <a< v_{1,2} <1,
\end{equation}
but when $a<0$, this no longer holds.\\
A bifurcation occurs when some of these fixed points coincide. We
have $v_1=v_2$ if and only if \begin{equation}\label{fn22b} \frac{(a-1)^2}{4} =\frac{b}{c},\end{equation}
and we have $v_0=v_1$ if \begin{equation}\label{fn22c} \frac{b}{c}=-a. \end{equation} When
$a=-1$ and $\frac{b}{c}=1$, all three fixed points coincide. The
case $a=-1$ is somewhat simpler because the cubic polynomial in
Eq.\eqref{fn20} reduces to \begin{equation}\label{fn22d} -v^3+v=\frac{b}{c}v, \end{equation} and
much of the subsequent analysis will only be carried out for that
particular case. For the moment, however, we return to the general
case. The linearization of Eq.\eqref{fn0a} at such a fixed point
$v_\star$ is
\begin{equation}\label{fn23}
\begin{split}
\left\{\begin{array}{lcl}
\varepsilon \frac{d\varphi}{d\tau}&=& -3v_\star^2\varphi +2(a+1)v_\star \varphi -a\varphi -\eta,\\\\
\frac{d\eta}{d\tau}&=& b\varphi -c \eta.
\end{array}\right.
\end{split}
\end{equation}
In order to determine the bifurcation behavior at such a fixed point, 
we need to investigate the eigenvalues of the Jacobian matrix given by
\begin{equation}\label{fn24}
J_{ij}=\left( \begin{array}{cc} \frac{1}{\varepsilon}(-3v_\star^2 +2(a+1)v_\star-a )&\:\:\:\:\: -\frac{1}{\varepsilon}\\\\
b&\:\:\:\:\: -c \end{array} \right).
\end{equation}
The stability of the fixed points $v_\star$ will depend on the signs of 
the trace and determinant of $J_{ij}$.
For a fixed point $v_\star$ to be stable, we should have $\mathrm{tr}J_{ij}<0$ and $\mathrm{det} 
J_{ij}>0$. Since $\varepsilon,c >0$, we have $\mathrm{tr}J_{ij}<0$ and $\mathrm{det} 
J_{ij}>0$ only if \begin{equation}\label{fn24a} -3v_\star^2 +2(a+1)v_\star  -a<0.
\end{equation} Eq.\eqref{fn24a} means that the fixed point $v_\star$ has to be
on the part of the cubic polynomial $-v^3+(a+1)v^2-av$ that
has a negative derivative (see the red curve in Fig.\ref{fig:chap32}\textbf{b}).
When we have three distinct fixed
points $v_\star$ in Eq.\eqref{fn21}, this could hold for the leftmost
and the rightmost of them,  and these two would then be  stable,
while the middle one would be unstable. The fixed point $v_0=0$ is
stable for $a>0$, and unstable for $a<0$.  We should point out,
however, that Eq.\eqref{fn24a} is a sufficient, but not a necessary
condition for the stability of a fixed point.

When we vary the parameters so that two of the fixed points merge,
we obtain a saddle-node type bifurcation. In the case $a=-1$, the
fixed points $v_{1,2}$ (which are symmetric in this case, that is,
$v_1=-v_2$) are stable according to Eq.\eqref{fn24a} if and only if 
\begin{equation}\label{fn24b}
v_1 < -\frac{1}{\sqrt{3}}, \end{equation} that is, if $v_1$  is to the left
of the local minimum of the cubic polynomial $w=-v^3+v$
(equivalently, $v_2>\frac{1}{\sqrt{3}}$ is to the right of the
local maximum). The limiting case \begin{equation}\label{fn24c} v_1
=-\frac{1}{\sqrt{3}}, \end{equation} of Eq.\eqref{fn24b} corresponds to
$\frac{b}{c}=\frac{2}{3}$, by Eq.\eqref{fn21}.  In this case, $v_1$
is precisely the minimum of the cubic polynomial $-v^3 +v$ in
Eq.\eqref{fn22d}. Likewise, $v_2$ then is the maximum of that
polynomial. This case will reoccur below in Eq.\eqref{fn29a}.

When Eq.\eqref{fn22} is not satisfied,  $v_0=0$ is the only fixed
point. For $v_0$, the determinant of Eq.\eqref{fn24} is
$\frac{1}{\varepsilon}(ac+b)$, and the trace is
$-\frac{a}{\varepsilon}-c$. Thus, $v_0$ is stable when $a$ is not
too negative, but in the limit $\varepsilon \to 0$, stability only
persists for $a\ge 0$. (Recall here that Eq.\eqref{fn24a} was
sufficient, but not necessary for the stability of a fixed point
$v_\star$.)

We obtain complex conjugate eigenvalues if
$0>\frac{1}{4}(\mathrm{tr} J_{ij} )^2 -\mathrm{det} J_{ij}$, and the real part of these eigenvalues vanishes when
$\mathrm{tr} J_{ij}=0$. For Eq.\eqref{fn24}, we thus have
the condition for complex conjugate eigenvalues as
\begin{equation}\label{fn25}
\begin{split}
\left\{\begin{array}{lcl}
\frac{1}{\varepsilon}(-3v_\star^2 +2(a+1)v_\star  -a) +c
-2\sqrt{\frac{b}{\varepsilon}}<0,\\
\frac{1}{\varepsilon}(-3v_\star^2 +2(a+1)v_\star  -a) +c
+2\sqrt{\frac{b}{\varepsilon}}>0,
\end{array}\right.
\end{split}
\end{equation}
and the condition for a
vanishing real part \begin{equation}\label{fn26} \frac{1}{\varepsilon}(-3v_\star^2
+2(a+1)v_\star  -a)-c=0. \end{equation} In order to be able to satisfy
Eq.\eqref{fn25} and Eq.\eqref{fn26} simultaneously  (and therefore have a 
Andronov-Hopf bifurcation), the coefficients
$b,c\ (>0)$ need to satisfy \begin{equation}\label{fn27} c^2 <\frac{b}{\varepsilon},
\end{equation} which is easily satisfied for small $\varepsilon >0$, since we
assume $b>0$.

When we solve Eq.\eqref{fn26} for $v_\star$, we obtain \begin{equation}\label{fn27aa}
v_{\star \star}=\frac{a+1}{3}\pm
\sqrt{\frac{(a+1)^2}{9}-\frac{a+\varepsilon c}{3}}, \end{equation} which we
can solve as long as \begin{equation}\label{fn27ab} 3\varepsilon c \le a^2-a +1. \end{equation}
We can then check when a solution $v_{\star \star}$ of
Eq.\eqref{fn27aa} coincides with one of the  points given by
Eq.\eqref{fn21}, in order to get a Andronov-Hopf bifurcation at one
of those  fixed points. We had already observed above that a fixed
point on the decreasing part of the cubic $v$-nullcline curve 
$w=-v^3+(a+1)v^2-av$ (the red curve in Fig.\ref{fig:chap32}\textbf{b})
is stable as the eigenvalues of the linearization then have
negative real parts. But this was a sufficient, but not a
necessary condition. Therefore, when we vary the slope
$\frac{b}{c}$ of the $w$-nullcline and move a fixed point to the middle increasing part
of the cubic polynomial, that fixed point may eventually lose its stability.
When the  parameter regime just identified for a Andronov-Hopf
bifurcation is right, that loss of stability will occur through a
Andronov-Hopf bifurcation.

In order to see the significance of Eq.\eqref{fn27ab}, we observe
that the local extrema of the cubic polynomial  $-v^3+(a+1)v^2-av$
are at \begin{equation}\label{rf27ac} v_\pm =\frac{a+1}{3}\pm
\sqrt{\frac{(a+1)^2}{9}-\frac{a}{3}}, \end{equation} and since $\varepsilon c
>0$, $\frac{a+1}{3}- \sqrt{\frac{(a+1)^2}{9}-\frac{a+\varepsilon
c}{3}}> \frac{a+1}{3}- \sqrt{\frac{(a+1)^2}{9}-\frac{a}{3}}$, that
is, the left Andronov-Hopf bifurcation is to the right of the
minimum of our cubic polynomial, hence on its ascending branch.
Whenever a fixed point is on the left descending branch, that is,
to the left of that minimum, it is stable, and stability persists
a little into the ascending branch, but when $\varepsilon \to 0$,
the Andronov-Hopf bifurcation point, that is, where the fixed point
looses its stability, moves towards that minimum. For $\varepsilon
>0$, the Andronov-Hopf bifurcation occurs to the right of the
minimum.

We return to the investigation of Eq.\eqref{fn26}. In particular, for
$v_\star =0$, Eq.\eqref{fn26} does not have a solution under our
constraint $ c>0$ when also $a>0$. When $a<0$, Eq.\eqref{fn26} is
satisfied for $v_\star =0$ when \begin{equation}\label{fn27a} a=-\varepsilon c.\end{equation}
When $a>0$ and Eq.\eqref{fn27ab} holds, the solutions $v_{\star
\star}$ of Eq.\eqref{fn27aa} are positive, and we might then tune the
parameter $b$ which does not occur in Eq.\eqref{fn27aa} so that one
of those $v_{\star \star}$ coincides with $v_1$ or $v_2$ from
Eq.\eqref{fn21}.

When we have equality in Eq.\eqref{fn27}, i.e., $c^2
=\frac{b}{\varepsilon}$, Eq.\eqref{fn27a} becomes \begin{equation}\label{fn27b}
a=-\frac{b}{c}, \end{equation} that is, Eq.\eqref{fn22c}. In that case, we have a
limit of a Andronov-Hopf bifurcation co-occurring with a
saddle-node bifurcation.

For $v_\star =v_{1,2}$ in Eq.\eqref{fn21}, we obtain 
\begin{eqnarray}\label{fn27b}\nonumber
2(a+1)\left(\frac{a+1}{2}\pm
\sqrt{\frac{(a-1)^2}{4} -\frac{b}{c}}\right)-a\\
=3\left(\frac{a+1}{2}\pm \sqrt{\frac{(a-1)^2}{4}
-\frac{b}{c}}\right)^2+\varepsilon c, 
\end{eqnarray}
 that is,
$$
 -\frac{a^2}{2} +a -\frac{1}{2} +3\frac{b}{c}-\varepsilon c=\pm (a+1)\sqrt{\frac{(a-1)^2}{4}-\frac{b}{c}},$$
which implies
\begin{eqnarray}\nonumber
\frac{1}{4}(a-1)^4 +(3\frac{b}{c} -\varepsilon c)^2 -(a-1)^2(3\frac{b}{c} -\varepsilon c)
\\=\frac{1}{4}(a-1)^2(a+1)^2 -(a+1)^2\frac{b}{c}\nonumber,
\end{eqnarray}
hence
\begin{eqnarray}\nonumber
-a^3+2a^2-a   -2a^2\frac{b}{c}+8a\frac{b}{c}-2\frac{b}{c} +a^2 \varepsilon c \\-2a\varepsilon c
+\varepsilon c +9\frac{b^2}{c^2} -6b\varepsilon +\varepsilon^2 c^2=0,\nonumber
\end{eqnarray}
which is equivalent to 
\begin{eqnarray}\label{fn28}\nonumber
(a-1)^2(-a
-2\frac{b}{c})+4a\frac{b}{c} +9\frac{b^2}{c^2}+\varepsilon^2 c^2\\=\varepsilon
(-(a-1)^2 c +6b).
\end{eqnarray}
%\end{equation} 
Thus, Eq.\eqref{fn28} is the
condition on the parameters for a Andronov-Hopf bifurcation at one
of the equilibria $v_{1,2}$.

In particular, in the limit $\varepsilon =0$, we get the condition
\begin{equation}\label{fn29} (a-1)^2(a +2\frac{b}{c})-4a\frac{b}{c}
=9\frac{b^2}{c^2}. \end{equation} For the case $a=-1$, we directly get from
Eq.\eqref{fn27b} \begin{equation}\label{fn29aa} 3\frac{b}{c}-2=\varepsilon c. \end{equation} Thus,
for $\varepsilon >0$, in Eq.\eqref{fn29aa}, we have $\frac{b}{c}
>\frac{2}{3}$, and therefore the fixed point is to the right of the
minimum, that is, on the increasing part of the cubic curve
$-v^3+v$, in accord with what we had said after Eq.\eqref{fn27aa}.

In the limit $\varepsilon =0$, we obtain \begin{equation}\label{fn29a}
\frac{b}{c}=\frac{2}{3}, \end{equation} as the condition for a Andronov-Hopf
bifurcation at an equilibrium $v_{1,2}$. This is the case of
Eq.\eqref{fn24c} (recalling Eq.\eqref{fn21}). Here, $v_1$ and $v_2$
lose their stability, and below in the slow-fast analysis, these
are also the points where the critical manifold will not be
normally hyperbolic, and where a switch from slow to fast dynamics
will occur, generating an unforced limit cycle. Since we had
identified Eq.\eqref{fn29a} as the parameter constellation for a
Andronov-Hopf bifurcation, this is a singular limiting situation
for a Andronov-Hopf bifurcation.

\begin{figure}%[H]
\begin{center}
\includegraphics[width=8.5cm,height=7.0cm]{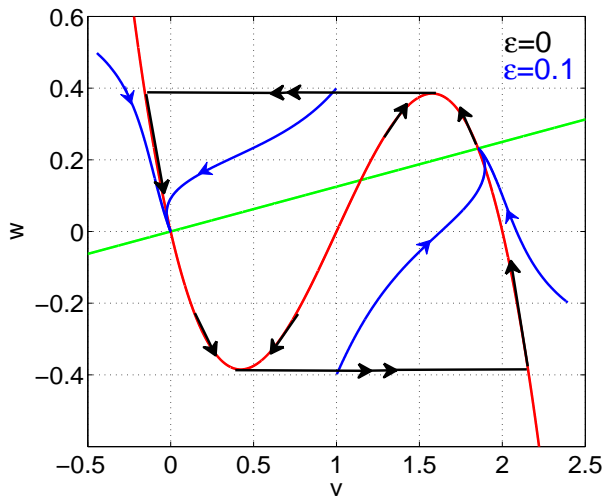}
\caption{(Color online) All the black trajectories (with single and double arrows) represent the solution of the slow flow of 
Eq.\eqref{fn44} (with $\varepsilon=0$) on the critical manifold $\mathcal{M}_0$ (the red curve).
The black trajectories with single arrow evolve on $\mathcal{M}_0$ towards 
the singular points at $v_\pm$ located at the extrema of $\mathcal{M}_0$, and horizontally leave
$\mathcal{M}_0$ at $v_\pm$ as shown by the black horizontal trajectories with double arrows.
The trajectories of the flow of Eq.\eqref{fn0a} (with $\varepsilon=0.1$) on the slow manifold $\mathcal{M}_{\varepsilon}$ (not shown but at a distance $\mathcal{O}(\varepsilon)$
from $\mathcal{M}_0$) converge to the stable fixed points $v_0=0$ and $v_2$ located at the intersection of the $w$-nullcline (the green line) and the decreasing parts of $\mathcal{M}_{0}$, 
without a limit cycle occurring.
$a=2.0, b=1.0, c=8.0, \sigma=0.0$} \label{fig:chap320}
\end{center}
\end{figure}

Now, we use slow-fast techniques to understand how a stable unforced 
limit cycle emerges from Eq.\eqref{fn0a} (with $\sigma=0$). 
For a more detailed introduction
to multiple time-scale dynamics see \cite{Kuehn}. In the singular limit 
$\varepsilon=0$, we define the critical manifold $\mathcal{M}_0$ of Eq.\eqref{fn0a} 
which coincides with our cubic polynomial nullcline: \begin{equation}\label{fn40}
\mathcal{M}_0:=\Big\{(v,w)\in\mathbb{R}^2:f(v,w)=0\Big\}. \end{equation} 
$\mathcal{M}_0$ can be viewed as the algebraic constraint of the differential-algebraic 
slow subsystem Eq.\eqref{fn44} whose initial conditions must satisfy
this constraint for solutions to exist. From Eq.\eqref{fn40}, we have \begin{equation}\label{fn41}
\frac{dw}{d\tau}= -3v^2 +2(a+1)v-a,
\end{equation}
with
\begin{eqnarray}\label{fn42}
\frac{dw}{d\tau}\begin{cases}<0,\,\, \text{for} \begin{cases}v<\frac{a+1}{3}-\frac{1}{3}\sqrt{a^2-a +1} \\
v> \frac{a+1}{3} +\frac{1}{3}\sqrt{a^2-a +1}
                \end{cases}\\ \\
=0,\,\,\text{for}\,\,v=\frac{a+1}{3} \pm \frac{1}{3}\sqrt{a^2-a +1} \\
                \\ 
>0,\,\, \text{for} \begin{cases} v>\frac{a+1}{3} -\frac{1}{3}\sqrt{a^2-a +1}\\
v< \frac{a+1}{3} +\frac{1}{3}\sqrt{a^2-a +1}
                \end{cases}
                \end{cases}
\end{eqnarray}
%\end{equation} 
Thus, $\mathcal{M}_0$ naturally splits into three parts: 
two decreasing stable branches and a middle
increasing unstable branch, see the red curve in Fig.\ref{fig:chap320} or
Fig.\ref{fig:chap32}\textbf{b}.  $\mathcal{M}_0$ 
looses its normal hyperbolicity at
the two singular points $v_\pm= \frac{a+1}{3} \pm \frac{1}{3}\sqrt{a^2-a+1}$, where it 
changes its stability property. These two points
are the local extrema of $\mathcal{M}_0$.
In fact, at $v_\pm= \frac{a+1}{3} \pm \frac{1}{3}\sqrt{a^2-a+1}$,
the existence and uniqueness theorems for ordinary differential 
equations (ODEs) do not longer apply, and because of this, the solutions of the slow subsystem 
in Eq.\eqref{fn44} are forced to leave $\mathcal{M}_0$
at these singular points.

In the singular limit $\varepsilon=0$, the dynamics of the fast variable $v$ on $\mathcal{M}_0$ (i.e., 
a 1-D dynamical system of the variable $v$ whose phase space is $\mathcal{M}_0$), is obtained from Eq.\eqref{fn40}  
by implicit differentiation \begin{equation}\label{fn43} (-3v^2
+2(a+1)v-a)\frac{dv}{d\tau}=\frac{dw}{d\tau}=bv-cw, \end{equation} and using the algebraic constraint on $w$ in Eq.\eqref{fn40}, 
we eliminate this variable to get the slow flow for the variable $v$ as
\begin{equation}\label{fn44} \frac{dv}{d\tau}=\frac{bv -c (-v^3+(a+1)v^2-av)}{-3v^2
+2(a+1)v-a}\ , \end{equation} which, of course, becomes singular at the points $v_\pm$.
In fact, at these points, since the critical manifold $\mathcal{M}_0$ loses its 
stability, the slow flow in Eq.\eqref{fn44} should detach from $\mathcal{M}_0$ and
should become fast, that is, satisfy $\frac{dw}{d\tau}=0$ and horizontally
jump to another stable branch of $\mathcal{M}_0$. 

In Fig.\ref{fig:chap320}, all the black trajectories (with single and double arrows) represent the solution of the slow flow of 
Eq.\eqref{fn44} on $\mathcal{M}_0$. Because of the failure of the existence and uniqueness theorems of ODEs at $v_-$,
the lower horizontal (fast) trajectory (in black with double arrow) 
leaves the left stable part of $\mathcal{M}_0$ at its local minimum at $v_-$ to 
the right stable part of $\mathcal{M}_0$. For the same reason, the upper horizontal (fast) trajectory (in black with double arrow) 
leaves the right stable part of $\mathcal{M}_0$ at its local maximum at $v_+$ to the left stable part of $\mathcal{M}_0$.
In this same figure, the non-horizontal (slow) trajectories (all in black with a single arrow) of Eq.\eqref{fn44}, 
evolve on $\mathcal{M}_0$ towards the singular points at $v_\pm$, where they eventually leave $\mathcal{M}_0$. 
 
The solution (or more precisely, the singular solution with $\varepsilon=0$) of the slow subsystem in Eq.\eqref{fn44} is related to the solution
of the full system Eq.\eqref{fn0a} with $\varepsilon>0$ by  Fenichel's theorem \cite{Feni1,Feni}.
By that  theorem, for $0<\varepsilon\ll1$, the  slow manifold $\mathcal{M}_{\varepsilon}$ is a perturbation of the critical manifold $\mathcal{M}_0$, and it has the following properties:
\begin{itemize}
 \item (F1) $\mathcal{M}_{\varepsilon}$ is diffeomorphic to $\mathcal{M}_0$.
 \item (F2) $\mathcal{M}_{\varepsilon}$ has  distance $\mathcal{O}(\varepsilon)$ from $\mathcal{M}_0$.
 \item (F3) The flow on $\mathcal{M}_{\varepsilon}$ converges to the slow flow on $\mathcal{M}_0$ as $\varepsilon \rightarrow 0$.
 \item (F4) $\mathcal{M}_{\varepsilon}$ is $C^r$-smooth for any $r<\infty$ (as long as $f,g\in C^\infty$).
 \item (F5) $\mathcal{M}_{\varepsilon}$ is normally hyperbolic, with the same stability properties w.r.t. the
fast variable $v$ as $\mathcal{M}_0$ (attracting, repelling or saddle-type).
 \item (F6) $\mathcal{M}_{\varepsilon}$ is usually not unique. Manifolds satisfying (F1)-(F5) lie at
distance $\mathcal{O}(e^-{K /\varepsilon})$ from each other $K>0$, $K=\mathcal{O}(1)$.
 \item (F7) Similar conclusions hold for the stable/unstable manifolds of $\mathcal{M}_0$.
\end{itemize}
Most importantly, the
flow of Eq.\eqref{fn0a} (with $0<\varepsilon\ll1$) on $\mathcal{M}_{\varepsilon}$  will 
follow the slow flow of Eq.\eqref{fn44} (with $\varepsilon=0$) on $\mathcal{M}_0$.
By the definitions in Eq.\eqref{fn1} and Eq.\eqref{fn40}, we see that the fixed points
$v_\star$ in 
Eq.\eqref{fn21} also lie on $\mathcal{M}_0$. We have the fixed points $v_0=0$ and $v_2$
 on the decreasing part of $\mathcal{M}_0$, and they are therefore stable, while the fixed point $v_1$ is 
located between $v_0=0$ and $v_2$ and  unstable. 

In Fig.\ref{fig:chap320}, the blue trajectories with arrows pointing in the direction of the flow
on the slow manifold $\mathcal{M}_{\varepsilon}$ (not shown, but at a distance $\mathcal{O}(\varepsilon)$ from $\mathcal{M}_0$), represent solutions
of Eq.\eqref{fn0a} with $\varepsilon=0.1$.
They converge towards to the stable fixed points $v_0$ and $v_2$ 
located respectively on the left and right stable decreasing parts of the slow manifold $\mathcal{M}_{\varepsilon}$.

For a better visualization,  we can  henceforth
discuss  the dynamics with respect to $\mathcal{M}_0$, since  Fenichel's theorem tells us that the same dynamical 
behavior takes place on $\mathcal{M}_{\varepsilon}$.
With the fixed point configuration: $v_0<v_1<v_2$ with $v_0$ and $v_2$ stable and  
located respectively on the left and right decreasing parts of  $\mathcal{M}_0$, $v_1$ unstable and located on the increasing part of $\mathcal{M}_0$,
trajectories of
Eq.\eqref{fn0a} cannot exhibit a spiking behavior (i.e., a limit cycle solution cannot emerge). This is because trajectories converge and stay at 
the stable fixed points $v_0$ and $v_2$ whenever they encounter them on 
decreasing parts of $\mathcal{M}_{0}$. Hence, these trajectories cannot evolve and reach the singular points $v_{\pm}$ located at the extrema of $\mathcal{M}_{0}$,
where they can jump from the left to the right and then from the right to the left part of  $\mathcal{M}_{0}$ 
to produce a limit cycle solution. In the blue trajectories in Fig.\ref{fig:chap320}, they stick at the fixed points $v_0=0$ and $v_2$.

Thus, it depends on the relative positions of the singular points at $v_\pm$ and the fixed
points at $v_0$ and $v_2$ on the critical manifold $\mathcal{M}_{0}$ whether the flow of Eq.\eqref{fn0a} will 
first reach a stable fixed point (and stay there with no possibility for a limit cycle) or 
first reach a singular point $v_{\pm}$ and detach from $\mathcal{M}_{0}$ with the emergence of a limit cycle.

When, say, $v_0$ is the smallest fixed point and is located to the
right of $v_-$ and $v_2$, the largest fixed point,  to the left of
$v_+$, we expect a periodic cycle for $0<\varepsilon\ll1$. For
instance, when we start close to the left stable branch of $\mathcal{M}_0$,
by Fenichel's theorem, the flow closely
and slowly follows $\mathcal{M}_0$ until we get into
the vicinity of $v_-$. There, $\mathcal{M}_0$ becomes unstable,
and the flow will move away from it and become fast and therefore
move not perfectly horizontally this time, but with some inclination (because $\varepsilon>0$) to the right until it comes into the
vicinity of the right stable branch of $\mathcal{M}_0$. It will
then again move slowly and close to $\mathcal{M}_0$ until it
gets into the vicinity of $v_+$. Hence it moves away fast, also not perfectly
horizontally but with some inclination to the left, until it encounters the left
stable branch again, and the cycle repeats.

We therefore see that 
even in the absence of a deterministic input current,
the slow-fast structure of Eq.\eqref{fn0a} (with $\sigma=0$) can
naturally induce a limit cycle solution (spiking) if
the flow in Eq.\eqref{fn0a} leaves from the neighborhood of the stable parts of  
$\mathcal{M}_0$ at the singular points $v_{\pm}$ before 
it encounters a stable fixed point. This can only
happen if the fixed point is located to right of $v_-$ or 
to the left of $v_+$.

Comparing Eq.\eqref{fn24a} and Eq.\eqref{fn41}, we see that the fixed
point $v_\star$ is stable for all $\varepsilon >0$ if it is to the
left of the minimum of $\mathcal{M}_0$. In that case,  
it is on the left stable branch of $\mathcal{M}_0$, and the 
dynamics on that stable branch will therefore converge towards $v_\star$, 
without a limit cycle emerging. Analogously, a fixed
point $v_\star$ is stable if it is on the right stable branch of 
$\mathcal{M}_0$ with no possibility for a limit cycle as well.
Again, however, these are sufficient, but
not necessary conditions for stability of the fixed points. 
For $\varepsilon>0$, stability persists a little into the middle 
(increasing) part of $\mathcal{M}_0$ as we have found when we discussed the
possibility of an Andronov-Hopf bifurcation. Therefore, to have a bi-stability 
regime consisting of a stable fixed point and a stable unforced limit cycle, 
the fixed point should be located on the middle unstable branch of $\mathcal{M}_0$ and 
has be to stable. This can be done by choosing the parameters
such that the fixed point is located between the minimum $v_-$ of $\mathcal{M}_0$ 
and its Andronov-Hopf bifurcation value. 

For the purpose of this work, we shall henceforth consider the
situation where $v_0=0$ is the only fixed point, that is we choose
$a$, $b$, and $c$ such that Eq.\eqref{fn22} is not satisfied. The
persistence of stability on the middle unstable part of $\mathcal{M}_0$
also occurs for $v_0=0$. For $a<0$, $v_0=0$ is to the right of the
minimum $v_-=\frac{a+1}{3} - \sqrt{\frac{(a+1)^2}{9}
-\frac{a}{3}}$ of $\mathcal{M}_0$ (i.e., $v_-<v_0$) and therefore in the region
where it eventually loses its stability through an Andronov-Hopf
bifurcation when $\varepsilon$ increases. It
now suffices to choose specific values of $a$ ($a<0$) and $c$ 
(both values of $a$ and $c$ also not satisfying Eq.\eqref{fn22}) such
that $v_0<-\frac{a}{c}=\varepsilon$ (see Eq.\eqref{fn27a}) to have a \textit{stable} fixed point
$v_0$ to the right of $v_-$.

With the stable fixed point $v_0=0$ on
the middle part of $\mathcal{M}_0$, and from the slow-fast
analysis above, we also have a stable limit cycle surrounding this fixed point.
For topological reasons these attractors should be separated from each 
other by a repeller, in this case an
unstable limit cycle. In fact, the unstable limit cycle is the boundary of the
basin of attraction of the stable fixed point $v_0$. This immediately indicates that
the Andronov-Hopf bifurcation of the fixed point $v_0$ which eventually occurs as
$\varepsilon$ increases should be sub-critical.

We choose and maintain throughout this work the values of the
parameters as: $a=-0.05$, $b=1.0$, and $c=2.0$. For these values, we have
$v_-=-0.25305<v_0=0<-\frac{a}{c}=0.025$ and therefore $v_0$ is located on the middle
part of $\mathcal{M}_0$ and it is stable. The Andronov-Hopf
bifurcation value of $\varepsilon$ is computed from Eq.\eqref{fn27a} as
$\varepsilon_{hp}=0.025$. 

For these values of the system
parameters, see Fig.\ref{fig:chap32}\textbf{a}, we computed the
bifurcation diagram by selecting the maximum values
of the action potentials $v$ as a function of the bifurcation
parameter $\varepsilon$. For $0.024\leq\varepsilon<0.025$, the fixed
point $v_0=0$ is unstable as $\mathrm{det}J_{ij}=\frac{0.9}{\varepsilon}>0$ and
$\mathrm{tr}J_{ij}=\frac{0.05}{\varepsilon}-2>0$ and surrounded by a stable
limit cycle, and therefore no bi-stability. At
$\varepsilon=\varepsilon_{hp}=0.025$, a sub-critical Andronov-Hopf bifurcation
occurs and the unstable fixed point $v_0=0$ changes its stability.
For $0.025<\varepsilon<0.027865$, the fixed point $v_0=0$ is stable ($\mathrm{det}J_{ij}>0$
and $\mathrm{tr}J_{ij}<0$ for those values of $\varepsilon$) and co-exist with
the stable limit cycle. Thus, for $0.025<\varepsilon<0.027865$ we have a
bi-stability regime. At $\varepsilon=\varepsilon_{sn}=0.027865$, we have a
saddle-node bifurcation of limit cycles, in which case the stable
limit cycle surrounding the stable fixed point $v_0=0$ shrinks and eventually 
collides with the boundary of the basin of attraction
of $v_0$ (i.e., the unstable limit cycle). In this saddle-node bifurcation, 
the stable and the unstable limit cycle
annihilate each other leaving behind the stable fixed point $v_0$.
The fixed point $v_0$ maintains its stability in the interval
$0.027865<\varepsilon\leq0.029$, within which we have no
bi-stability as there is only one attractor in the entire phase space.

In the bi-stability regime $\varepsilon_{hp}<\varepsilon<\varepsilon_{sn}$,
depending on whether the initial conditions are chosen in the basin of attraction of the fixed point or in that of the limit
cycle, the dynamics will converge to either the fixed point or to the
limit cycle. This behavior is seen in Fig.\ref{fig:chap32}\textbf{b} (also already seen in the time-series in Fig.\ref{fig:chap31}\textbf{a} and \textbf{b})
which shows a phase portrait of one trajectory with initial
conditions in the basin of attraction of the stable fixed
point at $(v_0,w_0)=(0,0)$ (the blue dot at the origin) and two other trajectories with 
initial conditions in the basin of attraction of the stable limit cycle
(the blue closed curve). In this work, we will focus on weak-noise effects on the 
spiking dynamics of Eq.\eqref{fn0a} with
$\varepsilon_{hp}<\varepsilon<\varepsilon_{sn}$.

%%%%%%%%%%%%%%%%%%%%%%%%%%%%%%%%%%%%%%%%%%%%%%%%%%%%%%%%%%%%%%%%%%%%%%%%%%%%%%%%%%%%%%%%%%%%%%%%%%%%%%%%%%%%
\section{Stochastic sensitivity analysis and the Mahalanobis metric}\label{sect4}
%%%%%%%%%%%%%%%%%%%%%%%%%%%%%%%%%%%%%%%%%%%%%%%%%%%%%%%%%%%%%%%%%%%%%%%%%%%%%%%%%%%%%%%%%%%%%%%%%%%%%%%%%%%%
\noindent
We now introduce noise to the neuron model (i.e., $\sigma>0$) and
perform a stochastic sensitivity analysis of the stable
attractors. In this section, we analyze the neuron model on the
fast time-scale $t$. The stochastic sensitivity matrix
associated to a stochastic dynamical system is an asymptotic
characteristic of the random attractors of the system
\cite{Bashkirtseva1}. For our model equation Eq.\eqref{fn0b} with
$0<\sigma\ll1$, it allows us to approximate a spread of random
trajectories around the stable fixed point $(v_0,w_0)=(0,0)$ and
stable limit cycle which we now denote by
$\big[\bar{v}(t),\bar{w}(t)\big]$. 
The random trajectories in the basin of attraction of the stable fixed point $(v_0,w_0)$ 
evolve according the evolution of the probability density of 
the FPE corresponding to Eq.\eqref{fn0b} \cite{Risken}. 

Suppose that a stationary solution, $P\big[v(t),w(t)\big]$, of this FPE exists. Generally, for
$n$-dimensional systems with $n\geq2$, one usually cannot  find such a stationary probability density
analytically \cite{Risken}. This is the situation with
Eq.\eqref{fn0b}. When $0<\sigma\ll1$, the constructive asymptotics
and approximations based on a quasi-potential function, $\varphi$,
given in Eq.\eqref{fn45} are frequently used \cite{M.I Freidlin}.
\begin{equation}\label{fn45} \varphi=-\lim_{\sigma \rightarrow 0}\sigma^2\log
P\Big\{\big[v(t),w(t)\big],\sigma\Big\}.\end{equation}

A quadratic form of the quasi-potential gives a Gaussian
approximation of $P_g\big[v(t),w(t)\big]$ in the vicinity of the
fixed point $(v_0,w_0)$,
\begin{align}\label{fn46}\nonumber
&P_g\Big\{\big[v(t),w(t)\big];(v_0,w_0)\Big\}\\&=\frac{1}{Z}
\exp\Bigg[-\frac{1}{2\sigma^2} \left(\begin{array}{c}
 v(t)-v_0\\w(t)-w_0 
\end{array}\right)^{\top}\Omega_{ij}^{-1}
\left(\begin{array}{c}
 v(t)-v_0\\w(t)-w_0 
\end{array}\right)\Bigg],
\end{align}
where $Z$ is the normalization constant and $\Omega_{ij}$ is
the covariance matrix of random trajectories around the stable fixed
point $(v_0,w_0)$, i.e., $\Omega_{ij}$ plays the role of the
stochastic sensitivity matrix of this stable fixed point and it is
determined by the algebraic equation \begin{equation}\label{fn47}
J_{ij}\Omega_{ij}+\Omega_{ij}J_{ij}^{\top}+G_{ij}=\textbf{0}, \end{equation}
$J_{ij}=\left( \begin{array}{cc} -a & -1\\
\varepsilon b & -\varepsilon c \end{array} \right)$ is the
Jacobian matrix at $(v_0,w_0)=(0,0)$ of the deterministic neuron
equation Eq.\eqref{fn0b}.
$G_{ij}=\left( \begin{array}{cc} 1 & 0\\
0 & 0 \end{array} \right)$ is the diffusion matrix of system
Eq.\eqref{fn0b} and $\top$ denotes the transpose.

As the fixed point $(v_0,w_0)$ is exponentially stable (that is, all
the eigenvalues of $J_{ij}$ at $(v_0,w_0)=(0,0)$ have strictly
negative real parts), the matrix equation in Eq.\eqref{fn47} has as its 
unique solution  the stochastic sensitivity matrix
$\Omega_{ij}\overline{}$ of the fixed point \cite{Bashkirtseva1}.
The eigenvalues $\lambda_k(\varepsilon)$, $k=\{1,2\}$ of the
stochastic sensitivity matrix $\Omega_{ij}$ define the variance
of the random trajectories around the fixed point $(v_0,w_0)$. The
largest eigenvalue (the largest SSF)
$\lambda_{max}=\text{max}\{\lambda_k(\varepsilon)\}$ for each
value of the singular parameter $\varepsilon$ indicates the
sensitivity of the stable fixed point $(v_0,w_0)$ to the random
perturbation. As $\lambda_{max}$ increases,
the sensitivity of the $(v_0,w_0)$ to noise also increases. This
means we have a higher probability of escape (i.e., shorter
residence time) from the basin of attraction of the stable fixed
point $(v_0,w_0)$ which we now denote for short as
$\mathcal{B}(v_0,w_0)$.

The matrix equation Eq.\eqref{fn47} for Eq.\eqref{fn0b} reduces to
the system of algebraic equations \begin{equation}\label{fn48}
\begin{split}
\left\{\begin{array}{lcl}
-2a\Omega_{11}-\Omega_{12}-\Omega_{21}+1=0,\\
b\varepsilon\Omega_{11}+(-a-c\varepsilon)\Omega_{12}-\Omega_{22}=0,\\
b\varepsilon\Omega_{11}+(-c\varepsilon-a)\Omega_{21}-\Omega_{22}=0,\\
b\varepsilon\Omega_{12}+b\varepsilon\Omega_{21}-2c\varepsilon\Omega_{22}=0.
\end{array}\right.
\end{split}
\end{equation}
In the parametric zone of bi-stability: $a=-0.05$, $b=1.0$,
$c=2.0$, $\varepsilon_{hp}<\varepsilon<\varepsilon_{sn}$, the
stochastic sensitivity matrix $\Omega_{ij}$ of the stable fixed
point at $(v_0,w_0)=(0,0)$ of Eq.\eqref{fn0b} is given by
\begin{equation}\label{fn49}
\Omega_{ij}=\Omega_{ji}=\left( \begin{array}{cc} \frac{4\varepsilon
+0.9}{3.6\varepsilon-0.09} & \frac{\varepsilon}{1.8\varepsilon-0.045}\\
\frac{\varepsilon}{1.8\varepsilon-0.045} & 
\frac{\varepsilon}{3.6\varepsilon-0.09} \end{array} \right),
\end{equation}
with the eigenvalues (the SSFs) given by
\begin{equation}\label{fn50a}
\lambda_{1,2}(\varepsilon)=\frac{5\varepsilon+
0.9\mp\sqrt{25\varepsilon^2+5.4\varepsilon +0.81}}{7.2\varepsilon-0.18},
\end{equation}
where $\lambda_{max}=\lambda_2(\varepsilon)$ .
The corresponding generalized eigenvectors are 
\begin{equation}\label{fn50b}
U_{1,2}(\varepsilon)=\left( \begin{array}{c} \frac{0.04\varepsilon}
{-0.3\varepsilon-0.9\mp\sqrt{25\varepsilon^2+5.4\varepsilon+0.81}} \\\\
0.05 \end{array} \right).
\end{equation}

\begin{figure}%[H]
\begin{center}
\textbf{(a)}\includegraphics[width=7.5cm,height=6.5cm]{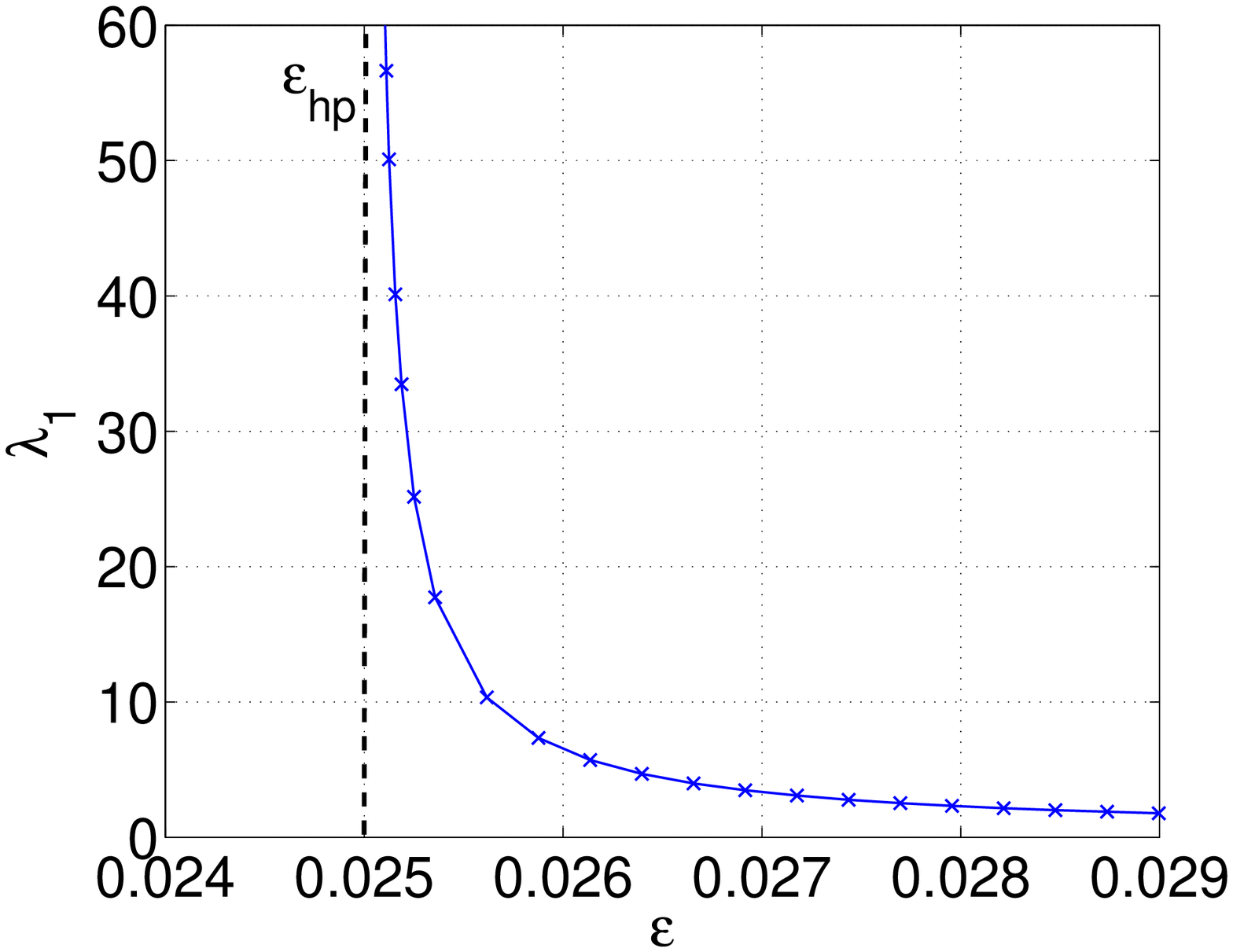}
\textbf{(b)}\includegraphics[width=7.5cm,height=6.5cm]{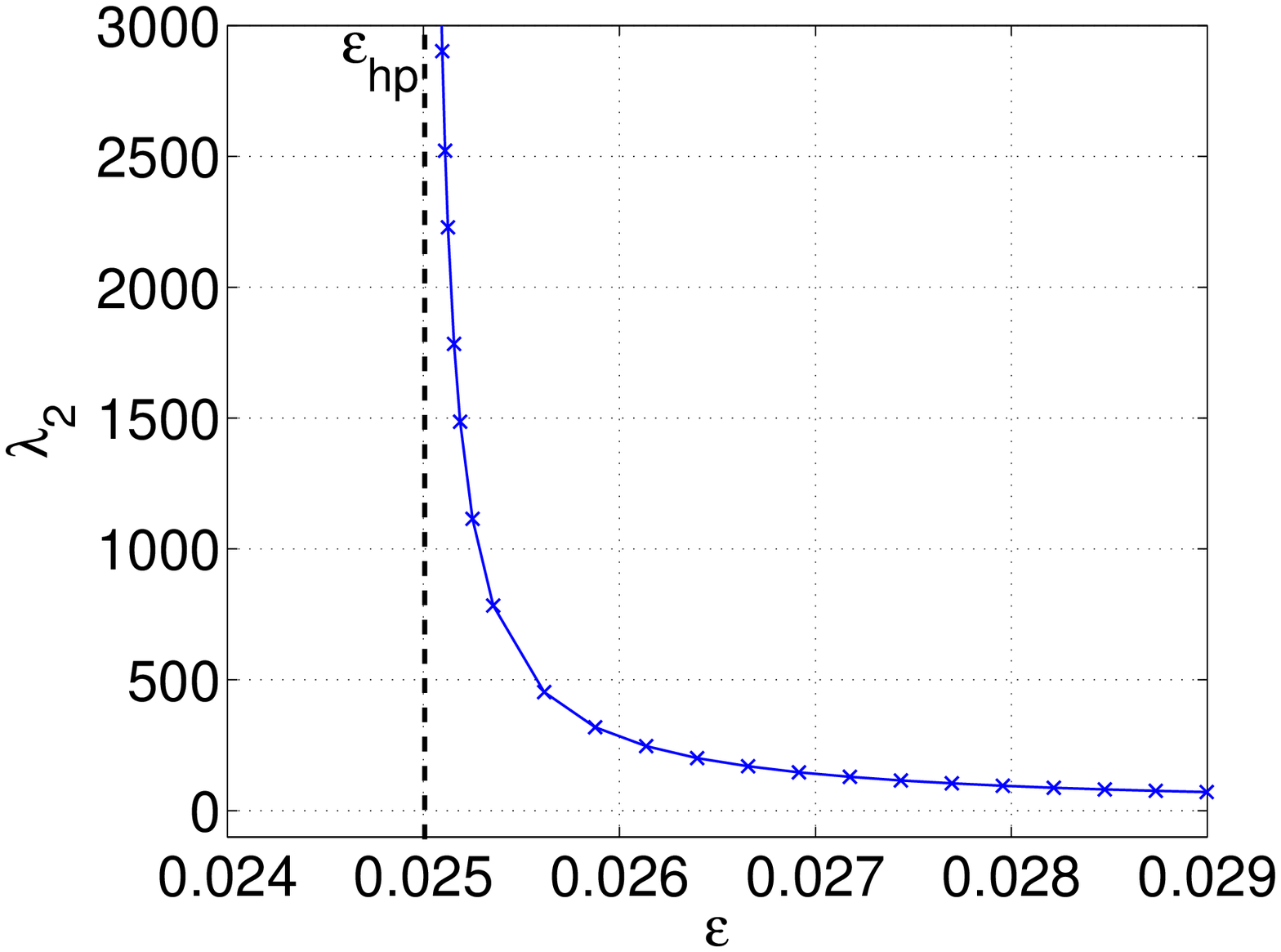}
\caption{Variation of the SSFs $\lambda_1$ in \textbf{(a)} and 
$\lambda_2$ in \textbf{(b)} of the stable fixed point at $(v_0,w_0)=(0,0)$ with 
the singular parameter $\varepsilon$. The SSFs diverge at
the Andronov-Hopf bifurcation value at
$\varepsilon=\varepsilon_{hp}=0.025$, with $\lambda_2$ 
dominating $\lambda_1$, indicating a higher stochastic sensitivity of the
fixed point in the direction of the corresponding
eigenvector $U_2$} \label{fig:chap33}
\end{center}
\end{figure}

For a fixed noise strength $\sigma$, the difference between
$\lambda_1$ and $\lambda_2$ reflects a spatial non-uniformity of
the dispersion of the random trajectories around the fixed point 
$(v_0,w_0)$ in the direction of the
eigenvectors $U_1$ and $U_2$ respectively. The dependence of
$\lambda_1$ and $\lambda_2$ on the singular parameter
$\varepsilon$ is shown in Fig.\ref{fig:chap33}.

Firstly, we
observe that the SSFs diverge as we approach the Andronov-Hopf
bifurcation value at $\varepsilon=\varepsilon_{hp}=0.025$. This means that the fixed point
$(v_0,w_0)$ becomes more and more sensitive to noise as we approach the
Andronov-Hopf bifurcation value and therefore the
highest probability of escape (i.e., shortest residence time) from $\mathcal{B}(v_0,w_0)$
when $\varepsilon\approx\varepsilon_{hp}$. 

Secondly, we observe that 
$\lambda_2$ diverges faster than $\lambda_1$ as $\varepsilon\rightarrow\varepsilon_{hp}=0.025$. This
shows that  the eigenvector $U_2$ localizes the main direction for
deviations of random trajectories from the fixed point
$(v_0,w_0)$, providing the direction in which the intersection
with the unstable limit cycle at $[v(t),w(t)]$ is most
probable. 

The Mahalanobis metric is a widely used metric in cluster and
discriminant analyses \cite{G. McLachlan}. Basically, it measures
the distance between a point $x$ and a distribution. This metric
is a natural tool for the quantitative analysis of noise-induced
transitions as it combines both the geometric distance from a
random attractor to a point and the stochastic sensitivity of this
attractor. The metric allows us to estimate a preference of the
stable fixed point $(v_0,w_0)$ or the stable limit cycle
$\big[\bar{v}(t),\bar{w}(t)\big]$ in the stochastic dynamics of
Eq.\eqref{fn0b}, when the random trajectory passes from one attractor
to another. For $0<\sigma\ll1$, the Mahalanobis distance from the
unstable limit cycle $\big[v(t),w(t)\big]$ (separatrix) to the
stable fixed point $(v_0,w_0)$ or to the stable limit cycle
$\big[\bar{v}(t),\bar{w}(t)\big]$ is related to the residence time
of trajectories in the corresponding basin of attraction: the
larger the Mahalanobis distance, the longer is the residence time
(i.e., lower probability of escape) in the corresponding basin.

In the stochastic sensitivity analysis of our fixed point
$(v_0,w_0)$, we approximate the probability density by a Gaussian
distribution in Eq.\eqref{fn46} centered at the stable fixed point at
$(v_0,w_0)=(0,0)$. The Mahalanobis distance
$D_m\Big\{\big[v(t),w(t)\big];(v_0,w_0)\Big\}$ from a point\\
$\big[v(t),w(t)\big]$ (i.e., a point on the unstable limit cycle) to
the distribution of random trajectories around the stable
attractor at $(v_0,w_0)$ is given by
\begin{align}\label{fn51}\nonumber
&D_m\Big\{\big[v(t),w(t)\big];(v_0,w_0)\Big\}\\
&=\sqrt{\left(\begin{array}{c}
 v(t)-v_0\\w(t)-w_0
\end{array}\right)^{\top}\Omega_{ij}^{-1}\left(\begin{array}{c}
 v(t)-v_0\\w(t)-w_0
\end{array}\right)},
\end{align}
where $\Omega_{ij}$ is the stochastic sensitivity matrix of
the fixed point $(v_0,w_0)$, and so the Gaussian approximation in
Eq.\eqref{fn46} can be written in terms of the
Mahalanobis distance as
\begin{align}\label{fn52}\nonumber
&P_g\Big\{\big[v(t),w(t)\big];(v_0,w_0)\Big\}\\  
&=\frac{1}{Z}\exp\left[-\frac{\bigg(D_m\Big\{\big[v(t),
w(t)\big];(v_0,w_0)\Big\}\bigg)^2}{2\sigma^2}\right].
\end{align}

For Eq.\eqref{fn0b}, we calculate the
Mahalanobis distances from the stable fixed point at $(v_0,w_0)=(0,0)$
to points on the unstable limit cycle at $\big[v(t),w(t)\big]$, and
then we choose the minimal distance. We calculate coordinates of
the unstable limit cycle numerically. This is done by assigning
$\big[v(t),w(t)\big]$ to the \textit{limiting values} of the initial
conditions $\big(v(0),w(0)\big)$ such that infinitesimal
perturbations (to the right and to the left) of these initial
conditions will lead to the convergence of the trajectories either
to the stable fixed point at $(v_0,w_0)$ or to the stable limit
cycle at $\big[\bar{v}(t),\bar{w}(t)\big]$ depending on which side
the infinitesimal perturbation is made.

We have
\begin{equation}\label{fn53}
\Omega_{ij}^{-1}=\left( \begin{array}{cc} 4\varepsilon-0.1 & -8\varepsilon+0.2\\
-8\varepsilon+0.2 & \frac{16\varepsilon^2+3.2\varepsilon-0.09}{\varepsilon} \end{array} \right),
\end{equation}
and using Eq.\eqref{fn51}, we calculate the
Mahalanobis distance from the stable fixed $(v_0,w_0)=(0,0)$ to
the unstable limit cycle $\big[v(t),w(t)\big]$ as
\begin{align}\label{fn54}\nonumber
&D_m\Big\{\big[v(t),w(t)\big];(0,0)\Big\}=\Bigg[ \frac{16\varepsilon^2+3.2\varepsilon-0.09}{\varepsilon}w(t)^2\\&
+(0.4-16\varepsilon)v(t)w(t)+(4\varepsilon-0.1)v(t)^2\Bigg]^{1/2}. 
\end{align}
Because of the dominance of $\lambda_2$ over $\lambda_1$, we
numerically calculate the minimum Mahalanobis distance $D_m$ from
the fixed point at $(v_0,w_0)=(0,0)$ to all points on the unstable
limit cycle at $\big[v(t),w(t)\big]$ by approximating the
Mahalanobis distance in Eq.\eqref{fn54} along the eigenvector
$U_2$, i.e.,
\begin{equation}\label{fn55}
D_m=\min_{(v,w)\in\big[v(t),w(t)\big]}D_m\Big\{\big[v(t),w(t)\big];(0,0)\Big\}.
\end{equation}

\begin{figure}[H]
\begin{center}
\includegraphics[width=7.5cm,height=6.0cm]{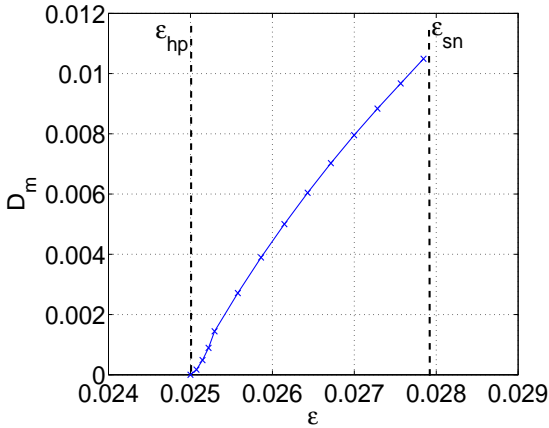}
\caption{Minimal Mahalanobis distance $D_m$ from the stable fixed
point at $(v_0,w_0)=(0,0)$ to the unstable limit cycle at
$\big[v(t),w(t)\big]$. Vertical dashed lines show the location
of the sub-critical Andronov-Hopf bifurcation of stable fixed point at
$\varepsilon_{hp}=0.025$ and of the saddle-node bifurcation of the stable
limit cycle at $\varepsilon_{sn}=0.027865$. $D_m$  vanishes at $\varepsilon=\varepsilon_{hp}$ and has maximum
value just before $\varepsilon=\varepsilon_{sn}$. $a=-0.05, b=1.0, c=2.0$} \label{fig:chap34}
\end{center}
\end{figure}

Fig.\ref{fig:chap34} shows the variation of the minimal
Mahalanobis distance $D_m$ from the fixed point $(v_0,w_0)$ to
the unstable limit cycle with the singular parameter $\varepsilon$.
The Mahalanobis distance vanishes at the sub-critical Andronov-Hopf
bifurcation value $\varepsilon_{hp}=0.025$, and increases with
increasing $\varepsilon$. This means as $\varepsilon$ increases
from $\varepsilon_{hp}$, the basin of attraction of the fixed
point increases in the direction of the eigenvector $U_2$, 
and therefore a lower and lower probability of
escape (i.e., longer residence time) from $\mathcal{B}(v_0,w_0)$ results.

From Figs. \ref{fig:chap33}\textbf{b} and \ref{fig:chap34}, we see that
the two factors (stochastic sensitivity of an  attractor and distance of the 
attractor to the separatrix) determining the length of the residence
time of random trajectories in $\mathcal{B}(v_0,w_0)$ are not 
competing. Approaching $\varepsilon_{hp}$ from above increases
the SSF of the fixed point and at the same time, decreases its
Mahalanobis distance to the unstable limit cycle. This has the combined
effect of considerably reducing the residence time of 
trajectories in $\mathcal{B}(v_0,w_0)$. In other words, there is a
higher probability (lower probability) that the random
trajectories escape from $\mathcal{B}(v_0,w_0)$ when $\varepsilon
\rightarrow \varepsilon_{hp}$ ($\varepsilon \rightarrow
\varepsilon_{sn}$).

We now apply a similar analysis to our randomly perturbed stable
limit cycle. In the deterministic system ($\sigma=0$) of
Eq.\eqref{fn0b}, with
$\varepsilon_{hp}<\varepsilon<\varepsilon_{sn}$, we have an
exponentially stable limit cycle defined by a $T$-periodic
solution,
$\big[\bar{v}(t),\bar{w}(t)\big]=\big[\bar{v}(t+T),\bar{w}(t+T)\big]$.
For the transversal hyperplane $\varSigma_{t}$ in the
neighborhood of any point $\big[\bar{v}(t),\bar{w}(t)\big]$ on the
stable limit cycle, the Gaussian approximation of the probability
density reads
\begin{align}\label{fn56}\nonumber
&P_g\Big\{\big[v(t),w(t)\big];[\bar{v}(t),\bar{w}(t)\big]\Big\}\\  \nonumber
&=\frac{1}{Z}\exp\Bigg[-\frac{1}{2\sigma^2}
\left(\begin{array}{c}
 v(t)-\bar{v}(t)\\w(t)-\bar{w}(t)
\end{array}\right)^{\top}\Theta_{ij}^{-1}(t)\\
&\times\left(\begin{array}{c}
 v(t)-\bar{v}(t)\\w(t)-\bar{w}(t)
\end{array}\right)\Bigg].
\end{align}

Here, the stochastic sensitivity matrix is periodic in time,
$\Theta_{ij}(t)=\Theta_{ij}(t+T)$. For an exponentially stable limit cycle, the
largest Lyapunov exponent is $0$ and the others are negative.
Consequently, the matrix $\Theta_{ij}(t)$ is the unique solution
of the Lyapunov equation \cite{Bashkirtseva1}, 
\begin{align}\label{fn57}\nonumber
 \frac{d\Theta_{ij}}{dt}&= J_{ij}(t)\Theta_{ij}(t)+\Theta_{ij}(t)J_{ij}(t)^{\top}+ P_{ij}(t) G_{ij} P_{ij}(t),
\end{align}
with the conditions
\begin{equation}\label{fn58}
\begin{split}
\left\{\begin{array}{lcl}
\Theta_{ij}(0)=\Theta_{ij}(T),\\\\
\Theta_{ij}(t)\left( \begin{array}{c} f\big[\bar{v}(t),\bar{w}(t)\big]\\\\
g\big[\bar{v}(t),\bar{w}(t)\big] \end{array} \right)\equiv 0,
\end{array}\right.
\end{split}
\end{equation} where $P_{ij}(t)$ is a matrix of the orthogonal projection
onto the Poincar\'{e} section $\varSigma_{t}$ at the point
$[\bar{v}(t),\bar{w}(t)]$ on the stable limit cycle, 
which is symmetric for our model equation Eq.\eqref{fn0b},
and whose entries are given by 
\begin{equation}\label{fn59}
\begin{split}
\left\{\begin{array}{lcl}
P_{11}=(-\bar{v}^3+(a+1)\bar{v}^2-a\bar{v}-\bar{w})^2,\\
P_{12}=\varepsilon(-\bar{v}^3+(a+1)\bar{v}^2-a\bar{v}-\bar{w})(b\bar{v}-c\bar{w}),\\
P_{22}=\varepsilon^2(b\bar{v}-c\bar{w})^2.
\end{array}\right.
\end{split}
\end{equation} 
$J_{ij}(t)$ is the Jacobian of the deterministic neuron at a point
$[\bar{v}(t),\bar{w}(t)]$ on the stable limit cycle and
given by \begin{equation}\label{fn60}
 J_{ij}(t)=\left( \begin{array}{cc} -3\bar{v}^2+2(a+1)\bar{v}-a &\:\:\:\:\: -1\\\\
\varepsilon b&\:\:\:\:\: -\varepsilon c\end{array} \right),
\end{equation}
and the constant diffusion matrix $G_{ij}$ is the same as before.

The Mahalanobis distance from a point 
$\big[v(t),w(t)\big]$ on the unstable limit cycle to the 
distribution of random trajectories around the stable
limit cycle at $\big[\bar{v}(t),\bar{w}(t)\big]$, 
$D_m\Big\{\big[v(t),w(t)\big];\big[\bar{v}(t),\bar{w}(t)\big]\Big\}$, is also a
periodic function of time and given by
\begin{align}\label{fn62}\nonumber
&D_m\Big\{\big[v(t),w(t)\big];\big[\bar{v}(t),\bar{w}(t)\big]\Big\}\\
&=\sqrt{\left(\begin{array}{c}
 v(t)-\bar{v}(t)\\w(t)-\bar{w}(t)
\end{array}\right)^{\top}\Theta_{ij}^{+}(t)\left(\begin{array}{c}
 v(t)-\bar{v}(t)\\w(t)-\bar{w}(t)
\end{array}\right)}.
\end{align}
Because $\Theta_{ij}(t)$ is singular for Eq.\eqref{fn0b}, 
``+'' means a pseudo-inverse in this case.

For 2-D systems, $\Theta_{ij}(t)$  can also
be written in the form \cite{Bashkirtseva2}
\begin{equation}\label{fn63b}
\begin{split}
\left\{\begin{array}{lcl}
\Theta_{ij}(t)=\mu(t)P_{ij}(t),\\
\Theta_{ij}^{+}(t)=\frac{1}{\mu(t)}P_{ij}(t),
\end{array}\right.
\end{split}
\end{equation}
and the Mahalanobis distance is given by 
\begin{equation}\label{fn63c}
 D_m\Big\{\big[v(t),w(t)\big];\big[\bar{v}(t),\bar{w}(t)\big]\Big\}=
 \frac{\left\lVert\left(\begin{array}{c}
 v(t)-\bar{v}(t)\\w(t)-\bar{w}(t)
\end{array}\right)\right\rVert}{\sqrt{\mu(t)}}.
\end{equation}
Here, $\mu(t)=\mu(t+T)>0$ is the unique solution of the boundary problem
\begin{equation}\label{fn64}
\begin{split}
\left\{\begin{array}{lcl}
d\mu=\alpha(t)\mu(t)dt+\beta(t)dt,\\
\mu(0)=\mu(T),
\end{array}\right.
\end{split}
\end{equation}
with $T$-periodic coefficients
\begin{equation}\label{fn64a}
\begin{split}
\left\{\begin{array}{lcl}
\alpha(t)=q(t)^\top\big[J(t)^\top + J(t)\big]q(t),\\
\beta(t)=q(t)^\top G_{ij}q(t).
\end{array}\right.
\end{split}
\end{equation}
$q(t)$ is a normalized vector orthogonal to the velocity vector field 
$\left( \begin{array}{c} f\big[\bar{v}(t),\bar{w}(t)\big] \\
g\big[\bar{v}(t),\bar{w}(t)\big]  \end{array} \right)$ and for
Eq.\eqref{fn0b} is given by
\begin{align}\label{fn65}\nonumber
&q(t)=\left(\begin{array}{c}-\varepsilon(b\bar{v}-c\bar{w})\\
-\bar{v}^3+(a+1)\bar{v}^2-a\bar{v}-\bar{w} \end{array} \right)\\
&\times\frac{1}{\sqrt{(-\bar{v}^3+(a+1)\bar{v}^2-a\bar{v}-\bar{w})^2
+\varepsilon^2(b\bar{v}-c\bar{w})^2}}.
\end{align}
We note that because our model is a $2D$ system, the hyperplane given by
$\varSigma_{t}$ is a tangent line to the stable limit cycle solution 
which is normal to $q(t)$ at
$[\bar{v}(t),\bar{w}(t)]$.
The functions $\alpha(t)$ and $\beta(t)$ for Eq.\eqref{fn0b} are given by
\begin{align} \label{fn66}\nonumber
&\alpha(t)=\frac{1}{(-\bar{v}^3+(a+1)\bar{v}^2-a\bar{v}-\bar{w})^2+\varepsilon^2(b\bar{v}-c\bar{w})^2}\\\nonumber
&\times\Bigg[2\varepsilon^2(-3\bar{v}^2+2(a+1)\bar{v}-a)(b\bar{v}-c\bar{w})^2\\\nonumber
&-2\varepsilon (\varepsilon b-1)(b\bar{v}-c\bar{w})(-\bar{v}^3+(a+1)\bar{v}^2-a\bar{v}-\bar{w})\\
&-2\varepsilon c(-\bar{v}^3+(a+1)\bar{v}^2-a\bar{v}-\bar{w})^2\Bigg].
\end{align}
\begin{align} \label{fn67}
\beta(t)=\frac{\varepsilon^2(b\bar{v}-c\bar{w})^2}{(-\bar{v}^3+(a+1)
\bar{v}^2-a\bar{v}-\bar{w})^2+\varepsilon^2(b\bar{v}-c\bar{w})^2}.
\end{align}

\begin{figure}%[H]
\begin{center}
\textbf{(a)}\includegraphics[width=7.5cm,height=5.7cm]{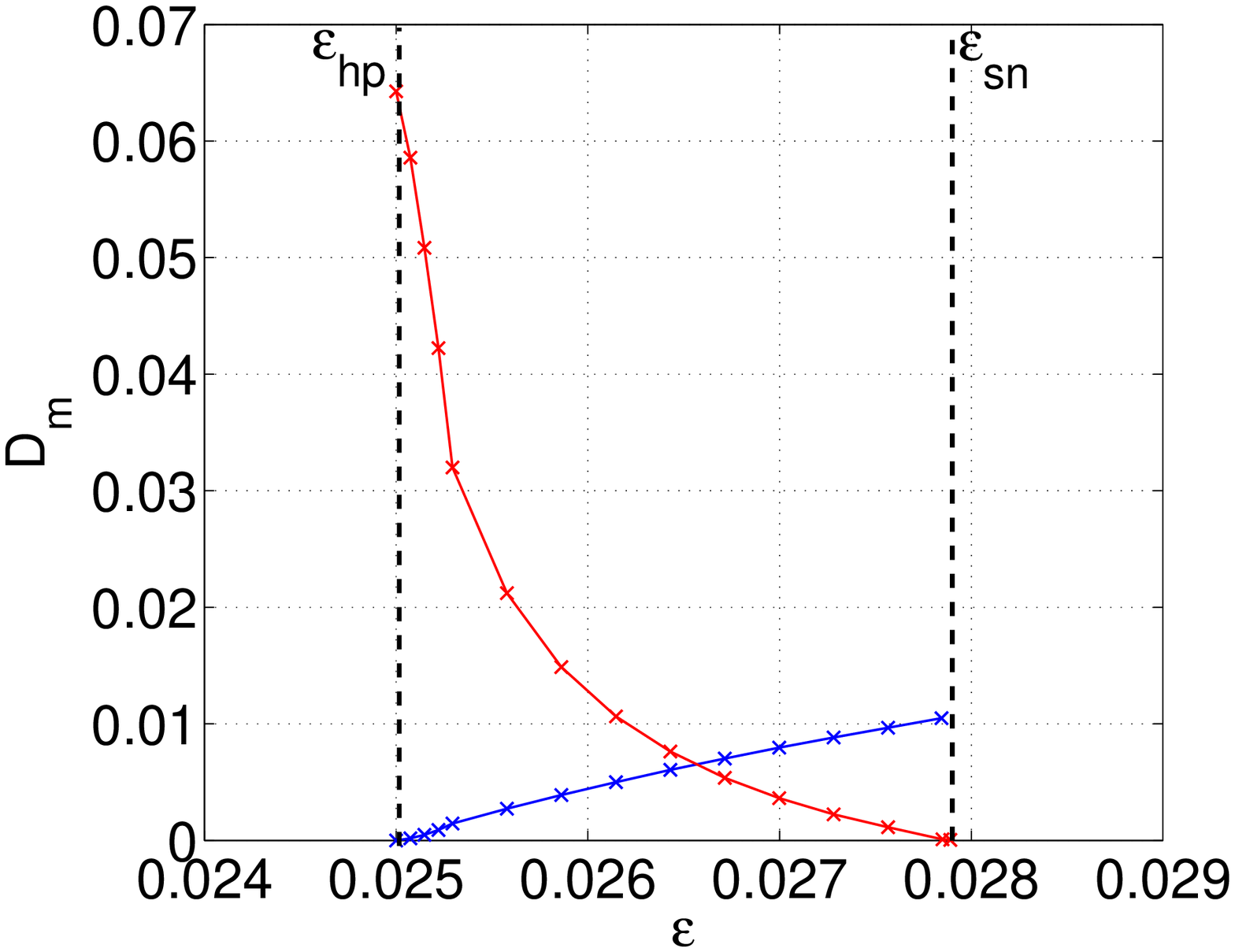}
\textbf{(b)}\includegraphics[width=7.5cm,height=5.7cm]{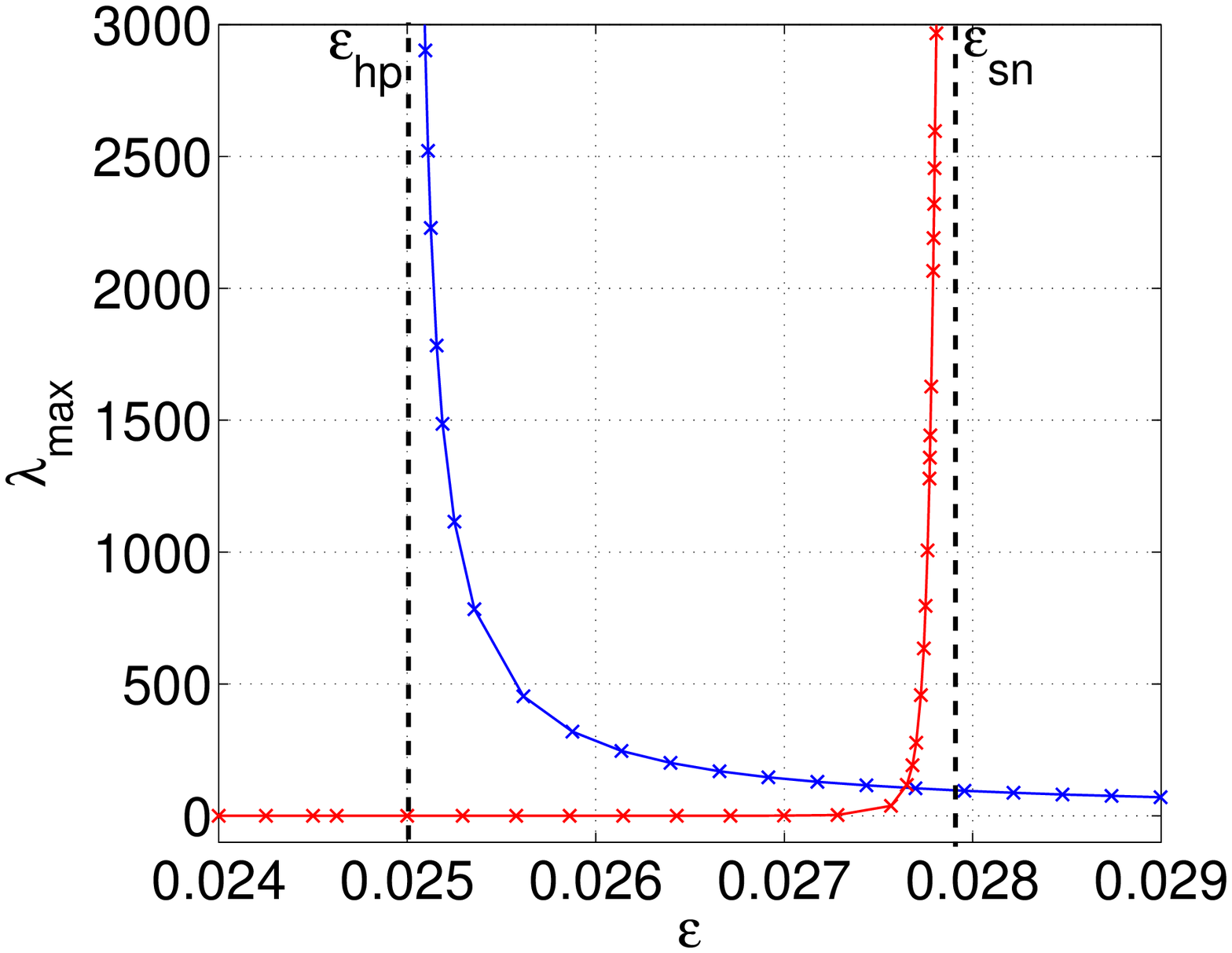}
\caption{(Color online) Variations of minimal Mahalanobis distances and the SSFs
of the attractors for Eq.\eqref{fn0b} with the singular parameter $\varepsilon$.
In \textbf{(a)}, the blue curve shows the minimal Mahalanobis distance, $D_m$, between stable
fixed point at $(v_0,w_0)=(0,0)$ and the unstable limit cycle at
$[v(t),w(t)]$ and the red curve corresponds to the minimal
Mahalanobis distance between stable limit cycle at
$\big[\bar{v}(t),\bar{w}(t)\big]$ and the unstable limit cycle, with $D_m$ vanishing at $\varepsilon_{sn}$.
\textbf{(b)} shows the SSFs of the stable fixed point (blue curve) and that of the stable limit
cycle (red curve).  
$a=-0.05, b=1.0, c=2.0$} \label{fig:chap35}
\end{center}
\end{figure}

The explicit solution of Eq.\eqref{fn64} is given by
\begin{equation}\label{fn68}
\mu(t)=e^{\int_0^t\alpha(s)ds}\Bigg[\int_0^t \beta(s)e^{\int_0^s-\alpha(r)dr}ds+C\Bigg].
\end{equation}
Because $\mu(t)$ is $T$-periodic, we write
\begin{align} \label{fn69}\nonumber
&e^{\int_0^t\alpha(s)ds}\Bigg[\int_0^t
\beta(s)e^{\int_0^s-\alpha(r)dr}ds+C\Bigg]\\\nonumber
&=e^{\int_0^{t+T}\alpha(s)ds}\Bigg[\int_0^{t+T}
\beta(s)e^{\int_0^s-\alpha(r)dr}ds+C\Bigg]\\ \nonumber
&=e^{\int_0^{t}\alpha(s)ds}e^{\int_t^{t+T}\alpha(s)ds}
\Bigg[\int_0^{t}\beta(s)e^{\int_0^s-\alpha(r)dr}ds\\
&+\int_t^{t+T}\beta(s)e^{\int_0^s-\alpha(r)dr}ds+C \Bigg],
\end{align}
and use the periodic property: $\alpha(t)=\alpha(t+T)$
with $\int_t^{t+T}\alpha(s)ds=\int_0^T\alpha(t+s)ds$, for
a fixed $t$, to get the constant $C$ as \begin{equation}\label{fn70}
C=\frac{e^{\int_0^{T}\alpha(s)ds}\cdotp
\int_0^{T}\beta(s)e^{\int_0^s-\alpha(r)dr}ds}{1-e^{\int_0^{T}\alpha(s)ds}}.
\end{equation}
With Eq.\eqref{fn63c} and the numerical value of $\mu(t)$ in Eq.\eqref{fn68}, 
the Mahalanobis distance is computed as
in Eq.\eqref{fn71} and the minimal Mahalanobis distance is calculated by taking
the minimum value of Eq.\eqref{fn71} over
$t\in[0,T)$ and $(v,w)\in[v(t),w(t)]$. See Fig.\ref{fig:chap35}\textbf{a}.
\begin{align}\label{fn71}\nonumber
D_m\Big\{&\big[v(t),w(t)\big];\big[\bar{v}(t),\bar{w}(t)\big]\Big\}\\
&=\sqrt{\frac{(v(t)-\bar{v}(t))^2+(w(t)-\bar{w}(t))^2}{\mu(t)}}.
\end{align}

This set, we obtain the
entries of $\Theta_{ij}(t)$ in Eq.\eqref{fn63b} using the
numerical value of $\mu(t)$ in Eq.\eqref{fn68}. As in the case
of the stable fixed point $(v_0,w_0)$, the eigenvalues
$\lambda_k(t)$, $k=\{1,2\}$, of $\Theta_{ij}(t)$ characterize
the distribution of random trajectories in the Poincar\'{e}
section $\varSigma_{t}$ near a point
$\big[\bar{v}(t),\bar{w}(t)\big]$ of the stable limit cycle. The
maximum of the largest eigenvalue indicates the SSF of the stable limit cycle. See Fig.\ref{fig:chap35}\textbf{b}.

%%%%%%%%%%%%%%%%%%%%%%%%%%%%%%%%%%%%%%%%%%%%%%%%%%%%%%%%%%%%%%%%%%%%%%%%%%%%%%%%%%%%%%%%%%%%%%%%%%%%%%%%
\section{Simulation results and discussion}\label{sect5}
%%%%%%%%%%%%%%%%%%%%%%%%%%%%%%%%%%%%%%%%%%%%%%%%%%%%%%%%%%%%%%%%%%%%%%%%%%%%%%%%%%%%%%%%%%%%%%%%%%%%%%%%
\noindent
In this section, numerical simulations are carried out with 
our model equation in the bistable regime to understand the
ISR we observed in Fig.\ref{fig:chap31}. We want to see how  ISR depends on which 
basin of attraction the initial conditions are located in, and how  the singular parameter
$\varepsilon$ affects ISR. We provide a theoretical explanation of the 
numerical results in terms of the results obtained in the stochastic
sensitivity analysis of our model equation. We recall 
that the differences in the SSFs and Mahalanobis distances
of our stable attractors define the direction of 
noise-induced transitions between them.

\begin{figure*}%[H]
\begin{center}
\textbf{(a)}\includegraphics[width=5.5cm,height=5.5cm]{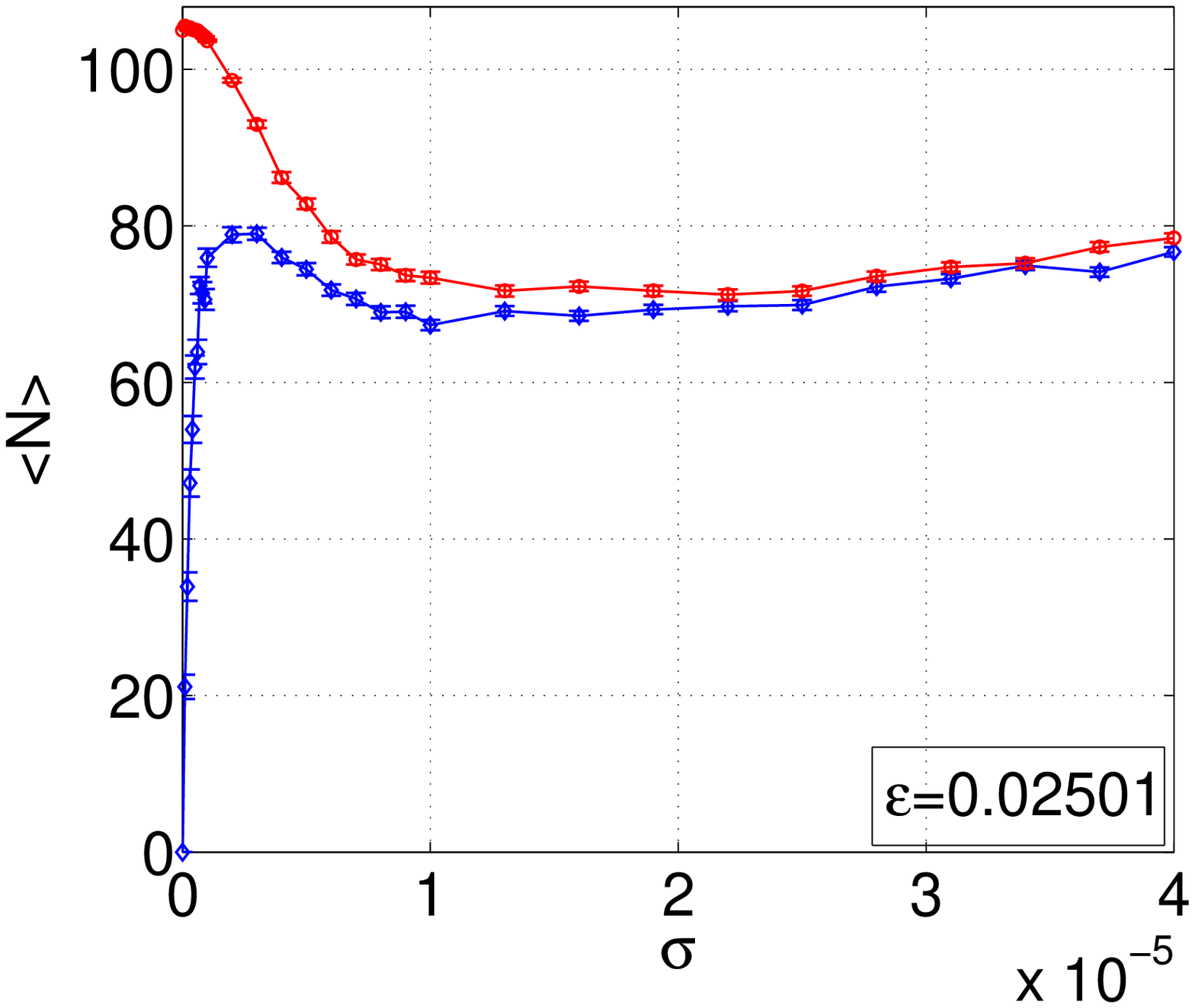}\textbf{(b)}\includegraphics[width=5.5cm,height=5.5cm]{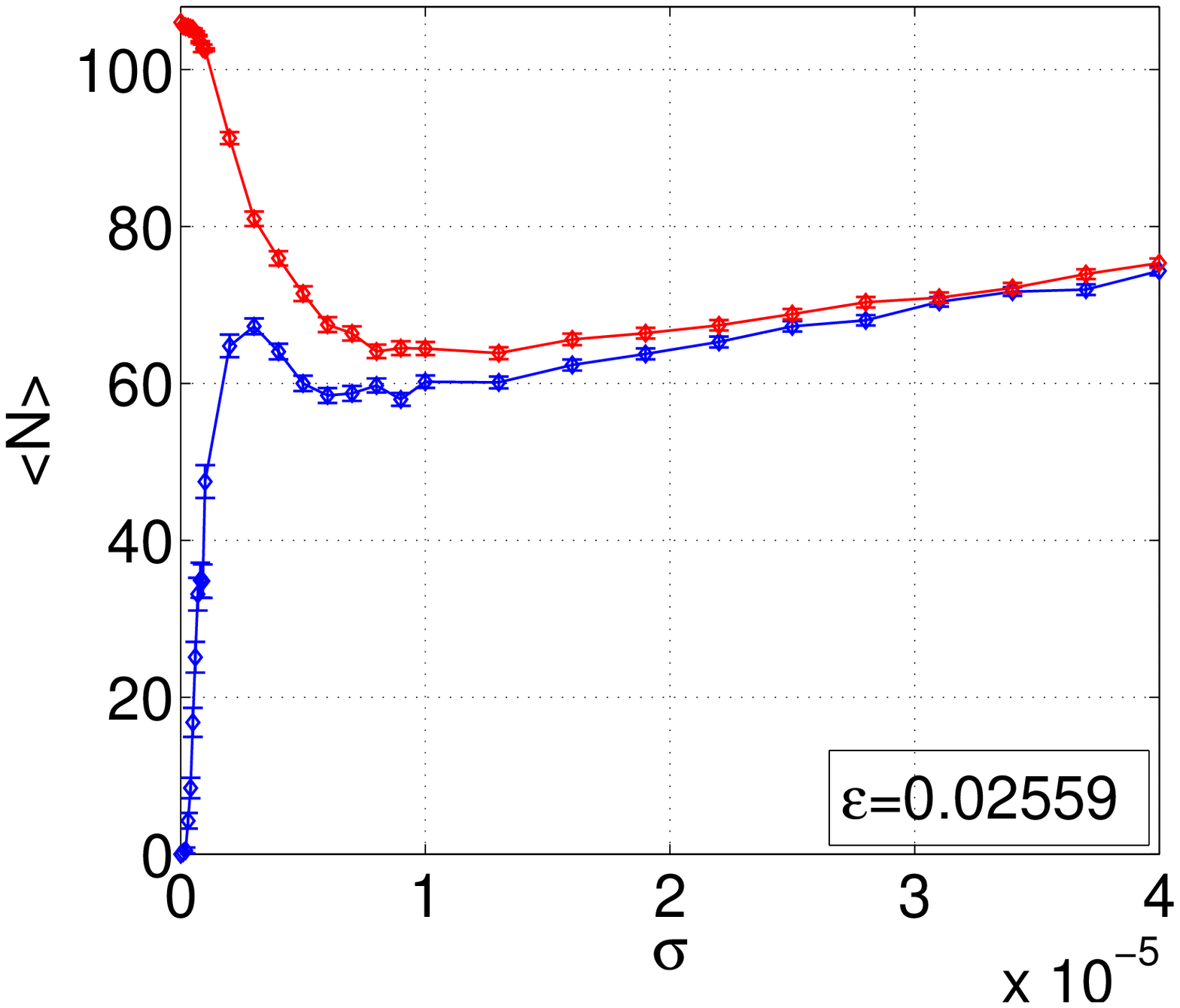}\textbf{(c)}\includegraphics[width=5.5cm,height=5.5cm]{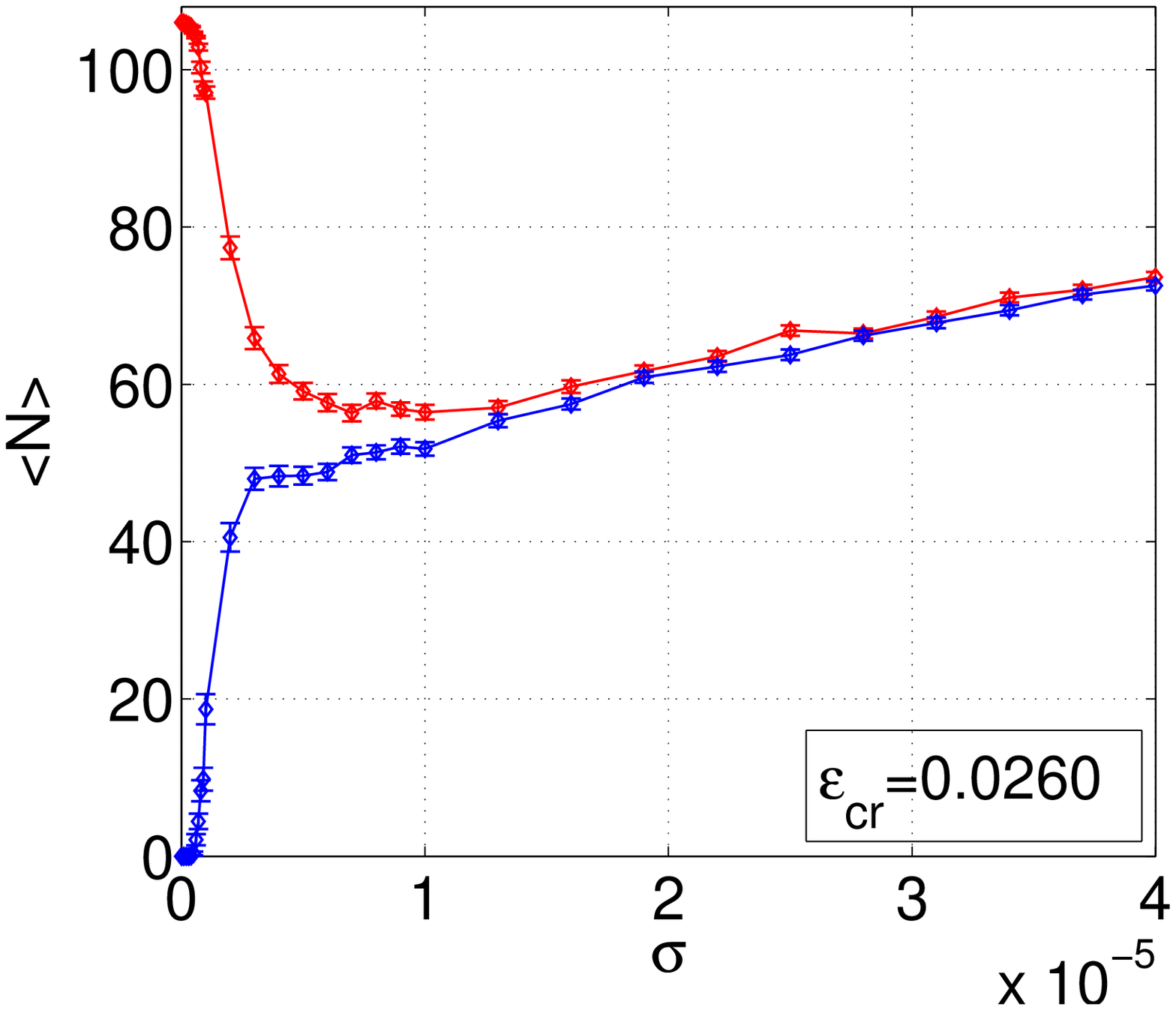}
\textbf{(d)}\includegraphics[width=5.5cm,height=5.5cm]{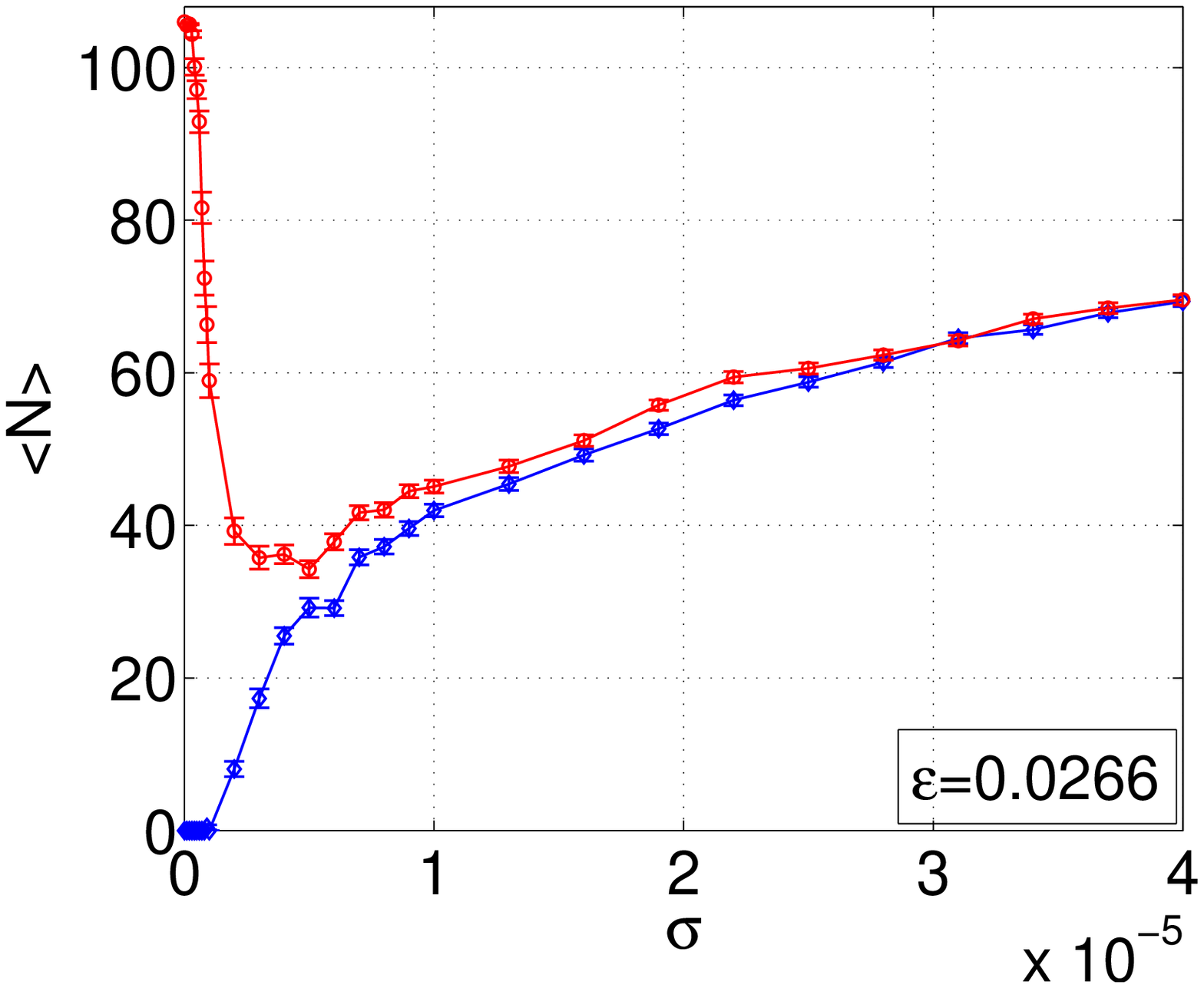}\textbf{(e)}\includegraphics[width=5.5cm,height=5.5cm]{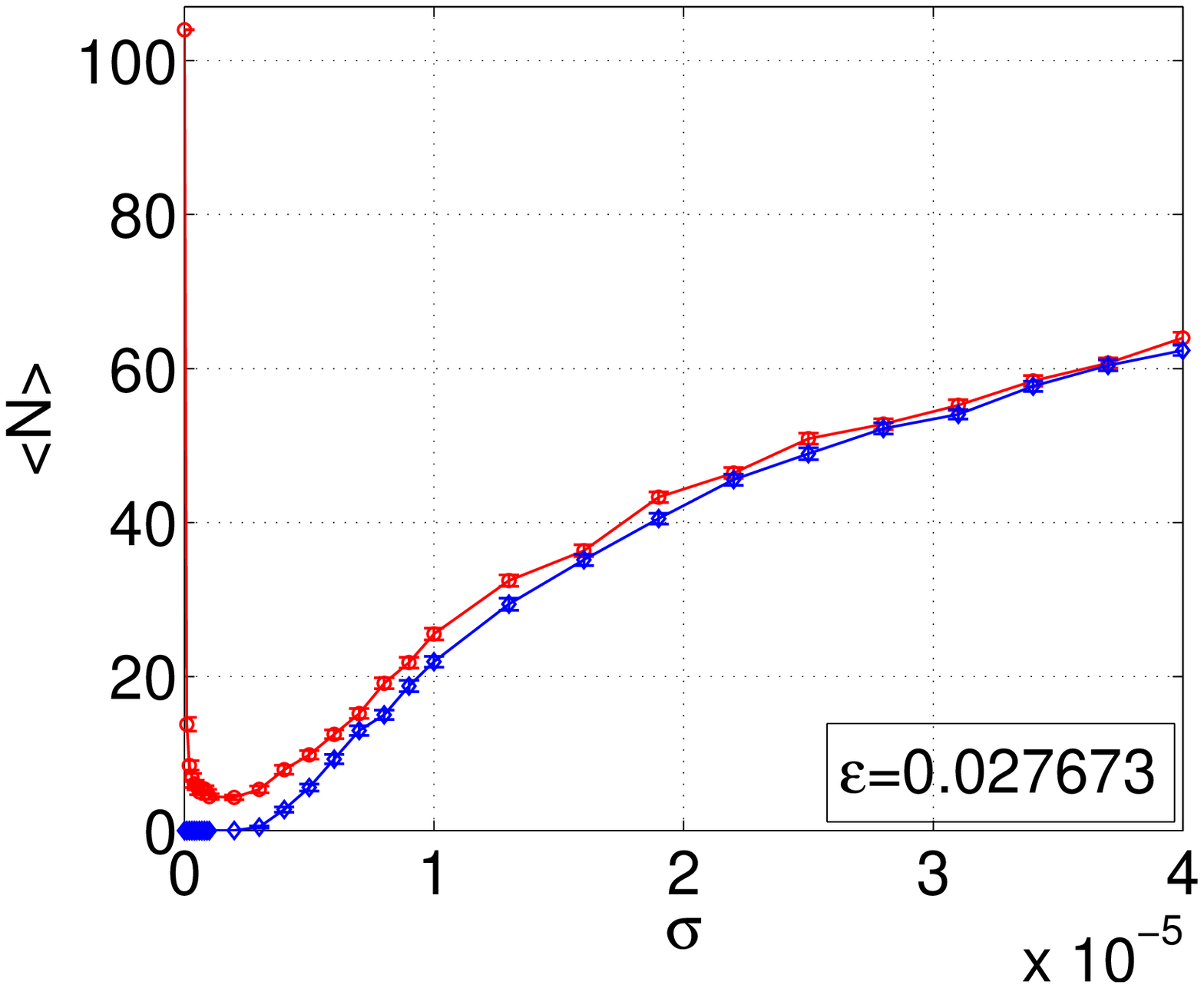}\textbf{(f)}\includegraphics[width=5.5cm,height=5.5cm]{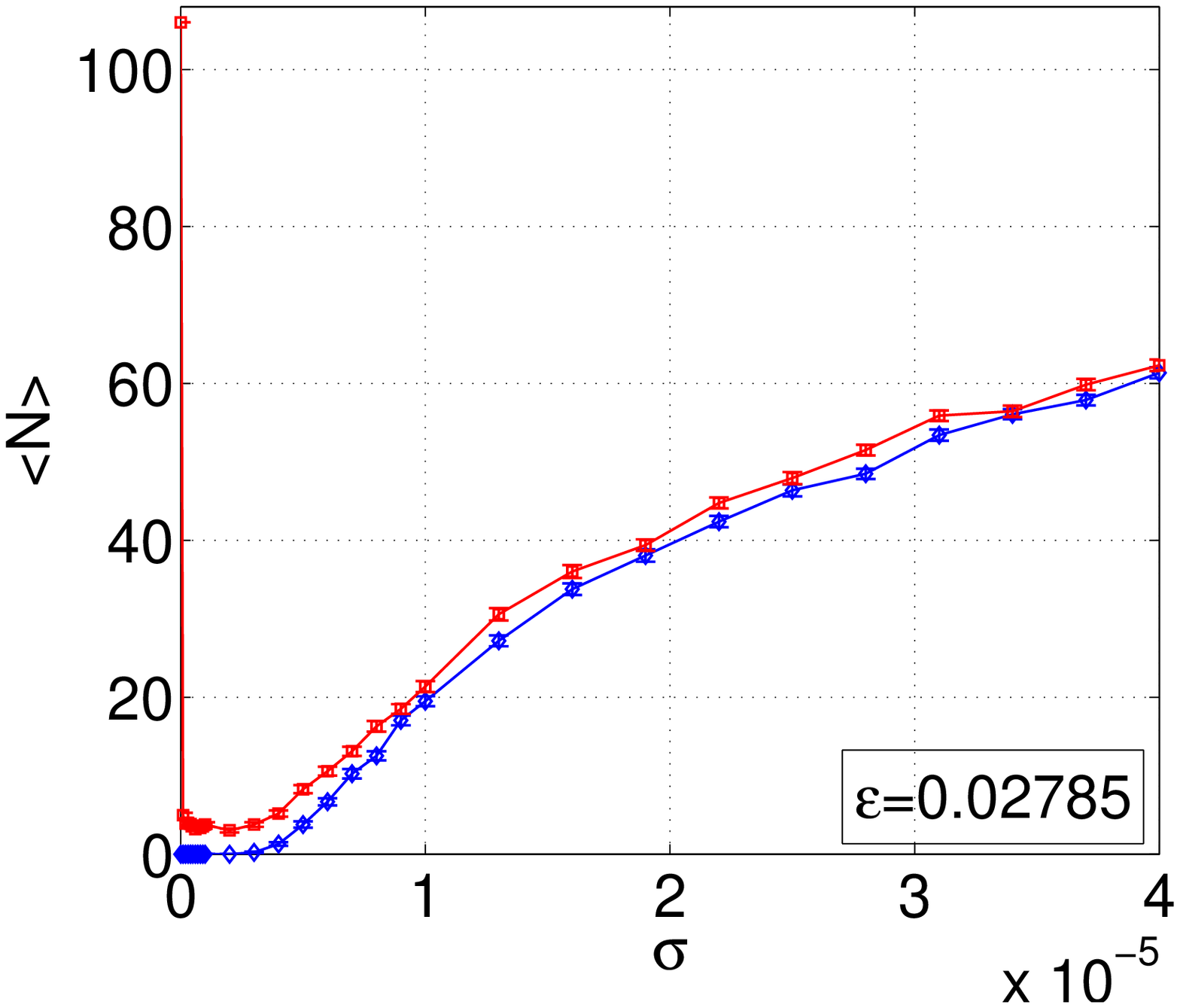}
\caption{(Color online) Mean number of spikes $\langle N\rangle$ versus noise
amplitude $\sigma$ ($200$ trials for $7500$ units of time interval each), for 
different singular parameter values $\varepsilon$ as indicated. ISR always occurs when
$\big(v(0),w(0)\big)\in \mathcal{B}\big[\bar{v}(t),\bar{w}(t)\big]$ (red curves).
For $\big(v(0),w(0)\big)\in \mathcal{B}(v_0,w_0)$ (blue curves), ISR \textit{only} occurs
when $D_m(fp)<D_m(lc)$, i.e, when $\varepsilon\in(0.025,0.0260)$.
$a=-0.05, b=1.0, c=2.0$. See text for details} \label{fig:chap36}
\end{center}
\end{figure*}

Using the fourth-order
Runge-Kutta algorithm for stochastic processes
\cite{Kasdin}, simulations are carried out for $200$
realizations of the noise and for $7500$ unit time intervals for 
each realization, a sufficiently long time interval for convergence of solutions
for $0<\sigma\ll1$. In Fig.\ref{fig:chap36}, 
we depict the variations of the mean number of
spikes $\langle N\rangle$ with the noise amplitude $\sigma$. The set of
numerical results are for different values of Mahalanobis
distances and SSFs of the stable attractors (encoded in the value of the singular parameter
$\varepsilon$ as in Fig.\ref{fig:chap35}). 

Sub-threshold responses are not counted as spikes. Again, we count a spike
(supra-threshold response) when the action potential
variable $v$ is greater than or equal to the threshold value of
$v_{th}=0.25$. We show simulation results for six values of the singular parameter
$\varepsilon\in(\varepsilon_{hp},\varepsilon_{sn})$, namely:
$\varepsilon=0.02501$ which is in the vicinity of the
Andronov-Hopf bifurcation value $\varepsilon_{hp}$ of the fixed
point, $\varepsilon=0.02559$, $\varepsilon=0.0260$,  
$\varepsilon=0.0266$ at which both attractors have equal
Mahalanobis distances, $\varepsilon=0.027673$
at which both attractors have equal SSFs, and
$\varepsilon=0.02785$ which is in the vicinity of the saddle-node
bifurcation value $\varepsilon_{sn}$ of the limit cycles.

The initial conditions $\big(v(0),w(0)\big)$ are fixed in every simulation. In Fig.\ref{fig:chap36}, the
red curves correspond to when $\big(v(0),w(0)\big)\in\mathcal{B}\big[\bar{v}(t),\bar{w}(t)\big]$. 
In this case, when $\sigma=0$, there
are $106$ spikes. The inhibitory effect of the spiking activity
begins as soon as $\sigma>0$, where we see the mean number of
spikes $\langle N\rangle$ decreasing to a minimum value before increasing
monotonically as $\sigma$ increases. We have ISR always occurring when 
$\big(v(0),w(0)\big)\in\mathcal{B}\big[\bar{v}(t),\bar{w}(t)\big]$.

The blue curves correspond to the situation where
$\big(v(0),w(0)\big)\in\mathcal{B}(v_0,w_0)$. In this case, when
$\sigma=0$, we have no spike, $\langle N\rangle=0$, and as soon as $\sigma>0$, 
$\langle N\rangle$ only increases monotonically and ISR does not occur. However, in
Fig.\ref{fig:chap36}\textbf{a} and \textbf{b} (blue curves), interestingly, ISR does actually occur
although $\big(v(0),w(0)\big)\in\mathcal{B}(v_0,w_0)$. 
We now  explain
these behaviors theoretically in terms of the Mahalanobis distances and the
SSFs of the attractors.

In Fig.\ref{fig:chap36}\textbf{a} (blue curve),
$\varepsilon=0.02501\gtrapprox\varepsilon_{hp}$, in which case the
Mahalanobis distance of the fixed point, $D_m(fp)$, is small and
smaller than the Mahalanobis distance of the limit cycle,
$D_m(lc)$. At this same value of $\varepsilon$, the SSF of the fixed point, $\lambda_{max}(fp)$, is
high and far higher than the SSF of
the limit cycle, $\lambda_{max}(lc)$. Therefore, with
$\big(v(0),w(0)\big)\in\mathcal{B}(v_0,w_0)$, very weak-noise
amplitudes are already capable of kicking the random trajectories out of
$\mathcal{B}(v_0,w_0)$ into
$\mathcal{B}\big[\bar{v}(t),\bar{w}(t)\big]$, thereby increasing
$\langle N\rangle$. For $0\leq\sigma\leq0.25\times10^{-5}$, $\langle N\rangle$ increases
from $0$ to a maximum of $79.8$. This same interval of $\sigma$ is
incapable of causing the reverse event, i.e., not strong enough to
kick random trajectories back into $\mathcal{B}(v_0,w_0)$ because of a
large $D_m(lc)$ and a low $\lambda_{max}(lc)$ and therefore, no
inhibitory effect of the spiking activity. As a
result, $\langle N\rangle$ can only  increase monotonically for
$0\leq\sigma\leq0.25\times10^{-5}$.
As soon as $\sigma>0.25\times10^{-5}$, it becomes strong enough to
kick the random trajectories it previously kicked into
$\mathcal{B}\big[\bar{v}(t),\bar{w}(t)\big]$ back to
$\mathcal{B}(v_0,w_0)$, thereby decreasing $\langle N\rangle$ (inhibiting the
spiking activity) down to a minimum value of $\langle N\rangle=69.3$ at
$\sigma=1.0\times10^{-5}$, before increasing monotonically with
$\sigma$.

From a series of simulations carried out for several different values of 
$\varepsilon\in(\varepsilon_{hp},\varepsilon_{cr}]$ (not all shown except for 
$\varepsilon=0.02501$, $\varepsilon=0.02559$, and $\varepsilon_{cr}=0.0260$ 
in Fig.\ref{fig:chap36}\textbf{a}-\textbf{c} respectively),
ISR remarkably persisted with $\big(v(0),w(0)\big)\in\mathcal{B}(v_0,w_0)$,
except at the critical value
$\varepsilon_{cr}$, where it just disappeared. 
That is, as $\varepsilon$ increased in the interval $(\varepsilon_{hp},\varepsilon_{cr})$ 
(with increasing $D_{m}(fp)$ and residence
time in $\mathcal{B}(v_0,w_0)$), ISR became less and less pronounced and eventually 
disappeared when $\varepsilon\geq\varepsilon_{cr}$.
This is so because as $D_{m}(fp)$ becomes larger and larger with increasing
$\varepsilon\geq\varepsilon_{cr}$, trajectories stay longer and longer in
$\mathcal{B}(v_0,w_0)$ (and therefore there is no spike), and as $\sigma$
increases, it becomes at each time \emph{just} strong enough to
kick the trajectories into $\mathcal{B}\big[\bar{v}(t),\bar{w}(t)\big]$
and hence $\langle N\rangle$ can only increase monotonically 
from $\langle N\rangle=0$ with $\sigma$. See
the blue curves in Fig.\ref{fig:chap36}\textbf{c}-\textbf{f} 
where $\varepsilon\geq\varepsilon_{cr}$, ISR does not occur as opposed to
the cases in Fig.\ref{fig:chap36}\textbf{a} and \textbf{b} (blue curves) where 
$\varepsilon\in(\varepsilon_{hp},\varepsilon_{cr})$.

Still in Fig.\ref{fig:chap36}\textbf{a} (now the red curve),
with $\big(v(0),w(0)\big)\\\in\mathcal{B}\big[\bar{v}(t),\bar{w}(t)\big]$, 
in the thin interval of very weak-noise amplitudes
$0\leq\sigma<0.2\times10^{-5}$, $\langle N\rangle$ is almost constant, near
$106$ (i.e., no considerable drop in $\langle N\rangle$ for this interval of
$\sigma$). This ``almost constant'' value of $\langle N\rangle$ in 
that interval of $\sigma$ happens because at
$\varepsilon=0.02501$, $D_{m}(lc)$ is large and
$\lambda_{max}(lc)$ is very low and as
$\big(v(0),w(0)\big)\in\mathcal{B}\big[\bar{v}(t),\bar{w}(t)\big]$,
trajectories have the tendency of staying in this basin of attraction for a very
long time and therefore almost no inhibition of the spiking
activity occurs for $0\leq\sigma<0.2\times10^{-5}$. As soon as
$\sigma>0.2\times10^{-5}$, it becomes strong enough to kick
trajectories out of the relatively larger
$\mathcal{B}\big[\bar{v}(t),\bar{w}(t)\big]$. The inhibitory effect
of noise then becomes pronounced with a clear decrease in $\langle N\rangle$
from about $106$ to a minimum of $78.3$ at
$\sigma=2.24\times10^{-5}$ before increasing monotonically with $\sigma$.

In Fig.\ref{fig:chap36}\textbf{d}, $\varepsilon=0.0266$, $D_{m}(fp)=D_{m}(lc)$, with
$\big(v(0),w(0)\big)\in\mathcal{B}\big[\bar{v}(t),\bar{w}(t)\big]$
(red curve), there is a rapid decrease in $\langle N\rangle$ for weak $\sigma$.
$\langle N\rangle$ moves from $106$ to a minimum of $34.6$ within
$0\leq\sigma\leq0.5\times10^{-5}$ before increasing monotonically
with increasing $\sigma$. A quicker decrease with a lower minimum
in $\langle N\rangle$ as compared to the cases in Fig.\ref{fig:chap36}\textbf{a}-\textbf{c} 
(red curves) occurs
because of a smaller $D_{m}(lc)$ in Fig.\ref{fig:chap36}\textbf{d} than in all 
previous cases, with
therefore a shorter residence time in
$\mathcal{B}\big[\bar{v}(t),\bar{w}(t)\big]$.

Still in Fig.\ref{fig:chap36}\textbf{d} with
$\big(v(0),w(0)\big)\in\mathcal{B}(v_0,w_0)$ (blue curve), ISR disappears.
Because at $\varepsilon=0.0266$ we have $D_{m}(fp)=D_{m}(lc)$, we explain this 
disappearance in terms of
the other factor determining the residence time in basins of
attraction i.e., the SSFs of the
attractors. At $\varepsilon=0.0266$, $\lambda_{max}(fp)$ is still
sufficiently high (with $\lambda_{max}(lc)<\lambda_{max}(fp)$,
which means that the fixed point is more sensitive to noise than
the limit cycle) and therefore, even weak-noise amplitudes have the
tendency of kicking the trajectories initially in
$\mathcal{B}(v_0,w_0)$ into
$\mathcal{B}\big[\bar{v}(t),\bar{w}(t)\big]$ (thereby increasing
$\langle N\rangle$). As in the cases of Fig.\ref{fig:chap36}\textbf{a} and \textbf{b} (blue curves),
one will expect that $\langle N\rangle$ increases with $\sigma\geq0$ up to a
certain maximum, and then start to decrease through the inhibitory
effect of noise. This is not happening in Fig.\ref{fig:chap36}\textbf{d} (and also already in 
Fig.\ref{fig:chap36}\textbf{c} blue curve)
firstly because $D_{m}(lc)$ is still sufficiently large (even though equal
to $D_{m}(fp)$) to keep the trajectories in
$\mathcal{B}\big[\bar{v}(t),\bar{w}(t)\big]$. Secondly, and mainly because 
$\lambda_{max}(lc)<\lambda_{max}(fp)$ at
$\varepsilon=0.0266$, which means that when trajectories get into $\mathcal{B}\big[\bar{v}(t),\bar{w}(t)\big]$,
they prefer to stay in this basin. And of course stronger
and stronger noise only increases $\langle N\rangle$.

In Fig.\ref{fig:chap36}\textbf{e}, $\varepsilon=0.027673$, $D_{m}(lc)$
is much smaller and $\lambda_{max}(lc)$ is much higher than in the
previous cases, but $\lambda_{max}(fp)=\lambda_{max}(lc)$. For the case
$\big(v(0),w(0)\big)\in\mathcal{B}\big[\bar{v}(t),\bar{w}(t)\big]$
(red curve), we therefore have a faster drop in $\langle N\rangle$, i.e., from
$106$ to a minimum of $6.4$ within
$0\leq\sigma\leq0.25\times10^{-5}$.  With
$\big(v(0),w(0)\big)\in\mathcal{B}(v_0,w_0)$ (blue curve) ISR
disappears for basically the same reason as previously given. In this
case, for $0\leq\sigma\leq0.25\times10^{-5}$, $\langle N\rangle$ remains at
zero (since $D_m(fp)>D_m(lc)$) and as $\sigma$ increases and becomes stronger, random trajectories start to jump into
$\mathcal{B}\big[\bar{v}(t),\bar{w}(t)\big]$ (thus increasing $\langle N\rangle$)
and remain in this basin for stronger and stronger noise with the
immediate consequence of just increasing $\langle N\rangle$.

In Fig.\ref{fig:chap36}\textbf{f},
$\varepsilon=0.02785\lessapprox\varepsilon_{sn}$, $D_{m}(lc)$ is
the smallest and $\lambda_{max}(lc)$ very high.
For $\big(v(0),w(0)\big)\in\mathcal{B}\big[\bar{v}(t),\bar{w}(t)\big]$
(red curve), there is a much more rapid and deeper drop in $\langle N\rangle$ with weak-noise
amplitudes as compared to all the previous cases with a well
defined minimum value of $\langle N\rangle=4.1$ at $\sigma=0.25\times10^{-5}$
and then a monotonic increase in $\langle N\rangle$ with increasing $\sigma$.
In this case, for some noise
realizations, the number of spikes could drop down to zero. That is, weak-noise amplitudes completely
terminate the spiking dynamics.

For $\big(v(0),w(0)\big)\in\mathcal{B}(v_0,w_0)$ (blue curve),
because $D_{m}(fp)$ is larger at $\varepsilon=0.02785$, the
residence time in $\mathcal{B}(v_0,w_0)$  is on average the
longest for weak-noise amplitudes compared to all previous cases. We have $\langle N\rangle=0$ for
$0\leq\sigma<0.30\times10^{-5}$. For $\sigma\geq0.30\times10^{-5}$, the noise is now
sufficiently strong to start kicking trajectories into $\mathcal{B}\big[\bar{v}(t),\bar{w}(t)\big]$
causing an increase in
$\langle N\rangle$. And as $\sigma$ becomes stronger, it keeps driving the neuron and so the
trajectories remain $\mathcal{B}\big[\bar{v}(t),\bar{w}(t)\big]$
with $\langle N\rangle$ increasing monotonically with $\sigma$.
%%%%%%%%%%%%%%%%%%%%%%%%%%%%%%%%%%%%%%%%%%%%%%%%%%%%%%%%%%%%%%%%%%%%%%%%%%%%%%%%%%%%%%%%%%%%%%%%%%%%%%%%%%
\section{Conclusion}\label{sect6}
%%%%%%%%%%%%%%%%%%%%%%%%%%%%%%%%%%%%%%%%%%%%%%%%%%%%%%%%%%%%%%%%%%%%%%%%%%%%%%%%%%%%%%%%%%%%%%%%%%%%%%%%%%%
\noindent
The effects of weak-noise amplitudes on the spiking dynamics of the FHN
neuron model without a deterministic input current was
investigated. Through bifurcation and slow-fast analyses, we
determined the conditions on the parameter space for the
establishment of a bi-stability regime consisting of a stable fixed
point and a stable unforced limit cycle. This bi-stability regime
induces a sensitivity to initial conditions in the
immediate neighborhood of the separatrix isolating the basins of attraction of the 
attractors. Introducing noise to the system then causes
transitions from the basin of attraction of the fixed point to that of
the limit cycle (the neuron gets into the spiking state) and as well, from the
basin of attraction of the limit cycle to that of the fixed point
(the neuron gets into the quiescent state, no spiking). 

We observed that in this bistable regime, weak-noise amplitudes
may decrease the mean number of spikes down to a minimum value
after which it increases monotonically as the noise strength increases. 
We showed that this phenomenon always occurred if the initial
conditions were chosen from the basin attraction of the stable
limit cycle. 

For initial conditions in the basin of attraction of
the stable fixed point, the phenomenon disappeared, unless  the time-scale separation parameter 
$\varepsilon$ of the neuron model
is bounded  between $\varepsilon_{hp}=0.0250$, 
its Andronov-Hopf bifurcation
value and $\varepsilon_{cr}=0.0260$. 
furthermore, the 
phenomenon became less and less pronounced
as $\varepsilon$ increased in the interval 
$(\varepsilon_{hp},\varepsilon_{cr})$, and
disappeared at $\varepsilon\geq\varepsilon_{cr}$.

We point out that this was not the case in \cite{Gutkin2} where
the neuron model considered had both a deterministic input and a
random input current. There, it was shown  that the decrease
to a minimum and then a monotonic increase in the spiking activity with
increase in noise amplitude occur regardless of the basin of
attraction from which the initial conditions are chosen 
from, provided the deterministic input current is above its
Andronov-Hopf bifurcation value. The model in the present work has
the same underlying dynamical structure as in \cite{Gutkin2}
except that we do not have a deterministic input current, and only  a
random perturbation component is considered.

We have seen that the stochastic sensitivity functions of the stable 
attractors and their Mahalanobis distances from the separatrix, which both
determine the length of the residence time of random trajectories in each 
state (quiescent or spiking state of the neuron),
themselves depended  on the time-scale 
separation parameter $\varepsilon$ of the model. 
From this dependence, we provided a theoretical 
explanation of the noise-induced phenomenon of ISR in
terms of the stochastic sensitivity functions and
the Mahalanobis distances of the stable attractors. 

Finally, we see that the key to ISR is the multi-stability between 
fixed points and limit cycles, a characteristic of dynamical systems with sub-critical
Andronov-Hopf bifurcations. To obtain bi-stability, in the present work, a careful relative 
positioning of the fixed point on the critical manifold was made. 
We can see in \cite{Yamakou1} how a small change in the relative position of fixed points
brings about a completely different dynamical behavior 
in the same weak-noise limit. Plausible implications of ISR in information processing
and transmission in neurons are discussed in \cite{Buchin}. 
\begin{acknowledgements}
This work was supported by the International Max Planck Research School 
Mathematics in the Sciences (IMPRS MiS),
Leipzig-Germany.
\end{acknowledgements}
%%%%%%%%%%%%%%%%%%%%%%%%%%%%%%%%%%%%%%%%%%%%%%%%%%%%%%%%%%%%%%%%%%%%%%%%%%%%%%%%%%%%%%%%%%%%%%%%%%%%%%%%%%%%%%%%%
%\section*{References}

\end{document}